\documentclass[a4paper,12pt]{article}
\usepackage[english]{babel}
\usepackage{enumitem}
\usepackage{color}
\usepackage{amsmath}
\usepackage{amssymb}
\usepackage{amsthm}
\usepackage{amsfonts}
\usepackage{graphicx}
\usepackage{epstopdf}
\usepackage{multirow}
\usepackage{cite}
\usepackage{caption}
\usepackage{graphicx} 
\usepackage{subcaption} 
\usepackage{amsmath, amssymb, mathtools, graphicx}
\usepackage{multirow}
\usepackage{array}
\usepackage{booktabs}
\usepackage{siunitx}
\usepackage{textcomp} 
\newcolumntype{Y}{>{\centering\arraybackslash}X} 

\newtheorem{theorem}{Theorem}[section]
\newtheorem{definition}{Definition}[section]
\newtheorem{lemma}{Lemma}[section]

\newtheorem{corollary}{Corollary}[section]
\newtheorem{remark}{Remark}[section]
\newtheorem{assumption}{Assumption}[section]
\newtheorem{algorithm}{Algorithm}[section]

\newcommand{\prox}{\operatorname{prox}}

\DeclareMathOperator*{\argmin}{argmin}

\DeclareMathOperator*{\dom}{dom}
\DeclareMathOperator*{\zer}{zer}

\DeclareMathOperator*{\Fix}{Fix}

\allowdisplaybreaks

\begin{document}
	
	\title{\textbf{ Preconditioned  Halpern iteration with adaptive anchoring parameters  and an  acceleration to Chambolle--Pock algorithm
	}}
	\author{ \sc \normalsize Fangbing Lv$^{a}${\thanks{ email: lfbmath@163.com}}\,\,\,and\,\,
		Qiao-Li Dong$^{a}${\thanks{Corresponding author.  email: dongql@lsec.cc.ac.cn}}\\
		\small $^a$College of Science,  Civil Aviation University of China, Tianjin 300300, China.\\
	}
	\date{}
	\maketitle
	
	\begin{abstract}
		In this article, we propose a preconditioned Halpern iteration with adaptive anchoring parameters (PHA) by integrating a preconditioner and Halpern iteration with adaptive anchoring parameters. Then we establish the strong convergence and at least $\mathcal{O}(1/k)$  convergence rate of the PHA method, and  extend these convergence results to Halpern-type preconditioned proximal point method with adaptive anchoring parameters. Moreover, we develop an accelerated Chambolle--Pock algorithm that is shown to have  at least $\mathcal{O}(1/k)$  convergence rate concerning the residual mapping and the primal-dual gap. Finally, numerical experiments on the minimax matrix game and LASSO problem are  provided to show the performance of our proposed algorithms.
	\end{abstract}
	
	
	\section{Introduction}
	Let \(U\) and \(V\) be  Hilbert spaces.  We consider the structured convex optimization problem formulated as
	\begin{equation}
		\label{tv1}
		\min_{u \in U}  f(u) + g(\mathbf{K}u),
	\end{equation}
	where  \(f: U \to \mathbb{R} \cup \{+\infty\}\) and \(g: V \to \mathbb{R} \cup \{+\infty\}\) are proper, convex and lower semi-continuous (lsc) functions. Here \(\mathbf{K}: U \to V\) is a bounded linear operator. The problem \eqref {tv1} has applications in numerous fields, including signal and image processing \cite{Bertsekas1982}, machine learning \cite{Bouwmans2016} and statistics \cite{Hayden2013}. Researchers introduced many methods to solve problem \eqref{tv1}, such as  Chambolle--Pock algorithm \cite{CP}, the alternating direction method of multipliers (ADMM) \cite{Gabay1976,Jonathan1992}, and their accelerated and generalized variants \cite{Liu2021,Malitsky}.
	
	The saddle-point or primal-dual form of  problem \eqref{tv1} reads as
	\begin{equation}\label{saddle}
		\min_{u \in U} \max_{v \in V} \mathcal{L}(u, v) := f(u) + \langle \mathbf{K}u, v \rangle - g^*(v),
	\end{equation}
	where $g^*$ is the Fenchel conjugate of $g$.
	The dual form of problem \eqref{tv1} is given by
	\begin{equation*}\label{fenchel}
		\max_{v \in V}  -g^*(v) - f^*(-\mathbf{K}^* v),
	\end{equation*}
	where \( \mathbf{K}^* \) is the adjoint operator of \( \mathbf{K}\). Under some regularity  conditions (see, e.g., \cite[Theorem 16.47]{Bauschke2017}),  the first-order optimality condition implies that the convex optimization problem (\ref{tv1}) is equivalently expressed as the inclusion problem
	\begin{equation}
		\label{tv2}
		0 \in \partial f(u) + \mathbf{K}^* \partial g(\mathbf{K}u),
	\end{equation}
	where \(\partial f\) and \(\partial g\) denote the subdifferentials of \(f\) and \(g\) respectively.  The problem \eqref{tv1}
	can also be reformulated as an one-operator inclusion problem
	\begin{equation}
		\label{monotone_inclusion}
		\text{find } x \in \mathcal{H} \text{ such that } 0 \in \mathcal{A}x,
	\end{equation}
	where $\mathcal{H}$ is a Hilbert space and $\mathcal{A}$ is maximal monotone. Indeed, letting \(\mathcal{H} = U \times V\), \(x = (u,v)\) and $
	\mathcal{A} = \begin{pmatrix}
		\partial f & \mathbf{K}^* \\
		-\mathbf{K} & \partial g^{*}
	\end{pmatrix}$,  \eqref{monotone_inclusion} is equivalent to \eqref{tv2}.
	
	By introducing a linear, bounded, self-adjoint, and positive semidefinite operator (so-called preconditioner)  $\mathcal{M} : \mathcal{H} \to \mathcal{H}$,  the inclusion problem \eqref{monotone_inclusion} can be equivalently recast as a fixed point problem
	\begin{equation}\label{4a1}
		x \in  (\mathcal{M} + \mathcal{A})^{-1} \mathcal{M}x.
	\end{equation}
	Assuming that the operator $(\mathcal{M} + \mathcal{A})^{-1}\mathcal{M}$ is single-valued and has full domain, one can solve \eqref{4a1}  by  Picard iteration
	\begin{equation}\label{4a}
		x^{k+1} = (\mathcal{M} + \mathcal{A})^{-1} \mathcal{M} x^k,
	\end{equation}
	which is called the degenerate preconditioned proximal point (PPP) method in \cite{Bredies}.
	The degenerate PPP method provides a unified framework and simplifies convergence analysis for some splitting methods used to solve monotone inclusion problems and structured convex optimization problems, such as the well-known Douglas--Rachford  algorithm \cite{Douglas1956}, Peaceman--Rachford algorithm \cite{Peaceman1955}, Davis--Yin method \cite{Davis2017},   proximal  ADMM algorithm \cite{Yang2025}  and the distributed forward-backward method  \cite{Aragon2023}. As noticed in \cite{Bredies}, Ryu's method \cite{Ryu2020}, Malitsky--Tam method \cite{Malitsky2023}, parallel Douglas--Rachford  algorithm \cite{Campoy2022} and  Chambolle--Pock  algorithm  can also be recast into the  degenerate PPP method \eqref{4a}.  Recently, the degenerate PPP method has given rise to new variants, such as the sequential Davis--Yin method \cite{Bredies}, the graph Douglas--Rachford \cite{Bredies2024} method and the graph forward-backward method \cite{Akerman2025,Aragón-Artacho2024}. Some authors have also investigated variable metric degenerate PPP methods \cite{Lorenz2024}, as well as inexact degenerate PPP methods with relative errors \cite{MarquesAlves2024}.
	
	Recently, the acceleration of the degenerate PPP method had drawn the attention of researchers.  Sun et al. \cite{Sun2025} proposed an accelerated degenerate proximal point method by integrating the  degenerate PPP method, Halpern iteration \cite{Halpern1967} and the fast Krasnosel'ski\v{\i}--Mann iteration \cite{Bot2024,Bot2023}, achieving asymptotic ${o}(1/k)$ and non-asymptotic $\mathcal{O}(1/k)$ convergence rates. In \cite{Zhang}, Zhang et al. introduced  a Halpern-type preconditioned proximal point (HPPP) algorithm which incorporates Halpern iteration with degenerate PPP method and has an \( \mathcal{O}(1/k) \) convergence rate. Moreover, they presented the Halpern-based accelerated Chambolle--Pock (HCP) method  by applying the HPPP algorithm to problem \eqref{tv1} with total variation (TV) regularization $\|\nabla u\|_1$ ($g=\|\cdot\|_1$, $\mathbf{K} = \nabla$)   as in \cite{Rudin1992}. However, the  convergence rate of the HCP algorithm  concerning the residual mapping or the primal-dual gap was not  provided.
	
	\vspace{2mm}
	\textbf{Our contributions}.
	Motivated by the above works, our contributions in this article are summarily given as
	\begin{enumerate}
		\item[{$\bullet$}]We propose a preconditioned Halpern iteration with adaptive anchoring parameters (PHA), which integrates a preconditioner into Halpern iteration with adaptive anchoring parameters \cite{He}. We analyze its strong convergence and give   at least $\mathcal{O}(1/k)$ convergence rate in Hilbert space. Furthermore, we extend the convergence results to Halpern-type preconditioned proximal point method with adaptive anchoring parameters.
		\item[{$\bullet$}]We develop an accelerated Chambolle--Pock (aCP) algorithm  by combining the PHA iteration and Chambolle--Pock algorithm. And we show that it achieves at least $\mathcal{O}(1/k)$ convergence rate concerning the residual mapping and the primal-dual value gap.
		\item[{$\bullet$}]We demonstrate the performance  of the proposed aCP algorithm and Restarted aCP algorithm through numerical experiments, including minimax matrix game and LASSO problem.
	\end{enumerate}
	
	\textbf{Organization of the article}.
	The rest of this article is organized as follows. Section 2 reviews some useful preliminaries for convergence analysis. Section 3 gives the PHA method and establishes its convergence results. In section 4, we propose  the aCP algorithm  and analyze its convergence rate. In section 5, numerical results are provided to show the performance  of the proposed algorithms. Finally, we conclude the article in Section 6.

	\section{Preliminaries}
	In this section, we start by introducing some concepts and lemmas for the convergence analysis.
	
	Let \( \mathcal{H} \) be a Hilbert space with the  inner product \( \langle \cdot, \cdot \rangle \) and the induced norm \( \|\cdot\| \). Denote by $I$  the identity operator on $\mathcal{H}$.  Let $\mathcal{M} : \mathcal{H}\rightarrow \mathcal{H}$ be a linear, bounded, self-adjoint and positive semidefinite operator.  For all $ x, y \in \mathcal{H}$, denote  the $\mathcal{M}$-induced semi-inner product and the $\mathcal{M}$-induced seminorm by $\langle x, y \rangle_{\mathcal{M}} = \langle \mathcal{M} x, y \rangle $ and  $ \|x \|_{\mathcal{M}} = \sqrt{\langle x, x  \rangle_{\mathcal{M}}} $, respectively.  We use the notations
	\begin{enumerate}
		\item[{(i)}] $\rightharpoonup$ for weak convergence and $\rightarrow$ for strong convergence.
		
		\item[{(ii)}] $\omega_w(x^k) = \{x : \exists\, x^{k_j}\rightharpoonup x\}$ denotes the weak $\omega$-limit set of $\{x^k\}$.
	\end{enumerate}
	
	\begin{definition}
		Let  $T: \mathcal{H}\rightarrow \mathcal{H}$ be an operator. Then $T$ is said to be
		\begin{enumerate}
			\item[{\rm(i)}] $\mathcal{M}$-nonexpansive if
			\[
			\|Tx - Ty\|_{\mathcal{M}} \leq \|x - y\|_{\mathcal{M}}, \quad \forall x, y \in \mathcal{H},
			\]
			and nonexpansive if $\|Tx - Ty\| \leq \|x - y\|,$  $\forall x, y \in \mathcal{H},$
			\item[{\rm(ii)}] $\mathcal{M}$-firmly nonexpansive if
			\[
			\|Tx - Ty\|_{\mathcal{M}}^2+\|(I-T)x-(I-T)y\|_{\mathcal{M}}^2 \leq \|x - y\|_{\mathcal{M}}^2, \quad \forall x, y \in \mathcal{H},
			\]
			and firmly nonexpansive if
			$
			\|Tx - Ty\|^2+\|(I-T)x-(I-T)y\|^2 \leq \|x - y\|^2, $ $\forall x, y \in \mathcal{H},$
			\item[{\rm(iii)}]
			$L$-Lipschitz continuous if
			\[
			\|Tx - Ty\| \leq L \|x - y\|,
			\]
and  $\mathcal{M}$-to-norm \(L\)-Lipschitz continuous if
			\[
			\|Tx - Ty\| \leq L \|x - y\|_{\mathcal{M}}.
			\]
		\end{enumerate}
	\end{definition}
	Denote by \( \Fix(T)\) the set of fixed points of \( T \), i.e., \( \Fix(T)= \{x \in \mathcal{H} : Tx = x\} \).

\begin{definition}
A set-valued mapping \(\mathcal{A}:\mathcal{H}\to 2^{\mathcal{H}}\) is said to be
\begin{enumerate}
			\item[{\rm(i)}] $\mathcal{M}$-monotone if
			\begin{equation*}
				\langle u - u', v - v' \rangle_{\mathcal{M}} \geq 0, \quad\forall   u \in \mathcal{A}v, \forall  u' \in \mathcal{A}v',
			\end{equation*}
			and monotone  if
			$
			\langle u - u', v - v' \rangle \geq 0, \quad\forall   u \in \mathcal{A}v, \forall  u' \in \mathcal{A}v',
			$
			\item[{\rm(ii)}] maximal monotone, if it is monotone and there exists no monotone operator whose graph  properly contains the graph of $\mathcal{A}$.
		\end{enumerate}
	\end{definition}
	Denote by $J_\mathcal{A} := (I +\mathcal{A} )^{-1}$ the resolvent of $\mathcal{A}$. Then $J_\mathcal{A}$  is firmly nonexpansive if $\mathcal{A}$ is maximal monotone. The set of solutions to (\ref{monotone_inclusion}) is referred to as zeros of \(\mathcal{A}\) and denoted by $\zer(\mathcal{A})$. By {\rm\cite[Proposition 20.36]{Bauschke2017}}, $\zer(\mathcal{A})$ is closed and convex  if $\mathcal{A}$ is maximal monotone.
	\begin{definition}
		Let $\mathcal{A}: \mathcal{H} \to 2^{\mathcal{H}}$ be a set-valued operator and $\mathcal{M}:\mathcal{H}\rightarrow \mathcal{H}$ be a  linear, bounded, self-adjoint and positive semidefinite operator.  $\mathcal{M}$ is said to be  an admissible preconditioner for $\mathcal{A}$ if
		$$
		(\mathcal{M} + \mathcal{A})^{-1} \mathcal{M} 
		\,\,\hbox{is single-valued and has full domain.}
		$$
	\end{definition}
	Suppose that $\mathcal{M}$ is an  admissible preconditioner. Then  according to {\rm\cite[Proposition 2.3]{Bredies}},  we have $\mathcal{M} = \mathcal{C} \mathcal{C}^*$, where $\mathcal{C} : \mathcal{D} \to \mathcal{H}$ is a bounded and injective operator for some real Hilbert space $\mathcal{D}$. Moreover, if $\mathcal{M}$ has closed range, then $\mathcal{C}^*$ is onto. Let $T=(\mathcal{M} + \mathcal{A})^{-1} \mathcal{M}$, then $\Fix(T) =\zer(\mathcal{A})$.
	
	Let $f: \mathcal{H} \to \mathbb{R} \cup \{+\infty\}$ be a proper, convex and lsc function. The subdifferential of $ f $ at $x$ is denoted by
	$$
	\partial f(x):=\left\{ u \in \mathcal{H} \mid  \langle y - x , u \rangle + f(x) \leq f(y),\,\,\forall y \in \mathcal{H} \right\}.
	$$
	Note that $\partial f$ is maximal monotone.
	The proximal operator of the function $ f $ with  $ \tau > 0 $ is   defined by
	$$
	\prox_{\tau f}(y) := \argmin_{x \in \mathcal{H}} \left\{ f(x) + \frac{1}{2\tau} \|x - y\|^2 \right\}, \quad \forall y \in \mathcal{H},
	$$
	which is firmly nonexpansive.
	Given \( x, p \in\mathcal{H}\), it holds
	\begin{equation}\label{d10}
		p = \prox_f(x) \iff x - p \in \partial f(p),
	\end{equation}
	which implies $\prox_f = J_{\partial f} $.
	The indicator function of a nonempty closed convex subset $D$ of $\mathcal{H}$ is denoted by $\iota_D$, i.e., $\iota_D(x) = 0$ if $ x \in D $ and  $\iota_D(x)=+\infty$ if $ x \notin D $. Note that $\prox_{\iota_D} = P_D$, where $P_D$ is the projector from $\mathcal{H}$ onto $D$.
	\begin{lemma}\label{lem211}
		Let $\mathcal{A} : \mathcal{H} \to 2^{\mathcal{H}}$ be a maximal monotone operator with $\zer(\mathcal{A}) \neq \emptyset$. Let $\mathcal{M}$ be an admissible preconditioner for $\mathcal{A}$ with a decomposition $\mathcal{M} = \mathcal{C} \mathcal{C}^*$ such that $(\mathcal{M} + \mathcal{A})^{-1}$ is $L$-Lipschitz continuous. Let $T=(\mathcal{M} + \mathcal{A})^{-1}\mathcal{M}$. Then
		\begin{enumerate}
			\item[{\rm(i)}] $T$ is $\mathcal{M}$-firmly nonexpansive;
			\item[{\rm(ii)}]$T$ is $\mathcal{M}$-to-norm \(L'\)-Lipschitz continuous with $L'=L \|\mathcal{C}\|$.
		\end{enumerate}
	\end{lemma}
	\begin{proof}
		(i) is derived by {\rm\cite[Proposition 2.5 and Lemma 2.6]{Bredies}}.
		(ii) Obviously, for $x, y \in \mathcal{H}$,
		\begin{equation*}
			\label{eq11}
			\begin{aligned}
				\|Tx - Ty\| &= \|(\mathcal{M} + \mathcal{A})^{-1}\mathcal{C}\mathcal{C}^*x - (\mathcal{M} + \mathcal{A})^{-1}\mathcal{C}\mathcal{C}^*y\|\\
				&\leq L \|\mathcal{C}\| \|x - y\|_{\mathcal{M}}.
			\end{aligned}
		\end{equation*}
		The proof is completed.
	\end{proof}
	
\begin{lemma}\label{n1}
Let $\mathcal{M} : \mathcal{H}\rightarrow \mathcal{H}$ be a linear, bounded, self-adjoint and positive semidefinite operator. Let $T$ be  $\mathcal{M}$-nonexpansive and  $\mathcal{M}$-to-norm \(L\)-Lipschitz continuous with a nonempty closed convex \( \Fix(T)\). Then, $\forall u \in \mathcal{H}$, there exists a unique point $u^*=\argmin_{v\in\Fix(T)}\|v-u\|_{\mathcal{M}}$ which solves
\begin{equation*}
\langle u - u^*, u^* - v \rangle_{\mathcal{M}} \geq 0, \quad \forall v \in \Fix(T).
\end{equation*}
\end{lemma}
\begin{proof}
We first establish the existence of  a minimizer for
    \begin{equation}\label{y1}
    \inf_{v \in \Fix(T)} \phi(v), \quad \phi(v) := \|v - u\|_{\mathcal{M}}.
    \end{equation}
Let \( \{v^k\} \subset \Fix(T)  \) be a minimizing sequence of $\phi$ ,  i.e., \( \phi(v^k) \to \inf_{v \in \Fix(T) } \phi(v) \) as $k \to \infty$. Fix some \( w \in \Fix(T)  \). Since \( T|_{\Fix(T)}  = \text{Id} \) and \( T \) is $\mathcal{M}$-to-norm \(L\)-Lipschitz continuous, we have
    \begin{equation}
			\label{eq111}
    \|v^k - w\| = \|T v^k - T w\| \leq L \|v^k - w\|_{\mathcal{M}}.
    \end{equation}
    Moreover,
    \begin{equation}
			\label{eq112}
    \|v^k - w\|_{\mathcal{M}} \leq \|v^k - u\|_{\mathcal{M}} + \|u - w\|_{\mathcal{M}}.
    \end{equation}
    Since \( \|v^k - u\|_{\mathcal{M}}= \phi(v^k) \) is bounded along the minimizing sequence, it follows from (\ref{eq112}) that \{\( \|v^k - w\|_\mathcal{M}\}\) is bounded, and hence from (\ref{eq111}) that  \( \{v^k\}\) is bounded. Therefore,  up to extraction of a subsequence, \( v^{k_l} \rightharpoonup v^* \) as $l \to \infty$ . As \( \Fix(T) \) is closed and convex, it is weakly closed, so \( v^* \in \Fix(T) \). Finally, the function \( \phi\) is convex and continuous, hence weakly lower semicontinuous. Consequently
    \[
    \phi(v^*) \leq \liminf_{l \to \infty} \phi(v^{k_l}) = \inf_{v \in \Fix(T)} \phi(v),
    \]
    and \( v^* \) attains the minimum.

Next we show the uniqueness of  the minimizer for (\ref{y1}). It is easy to see that
$u^*=\argmin_{v\in\Fix(T)}\|v-u\|_{\mathcal{M}}$ is equivalent  to
$$
u^*=\argmin_{v \in \mathcal{H}} \frac12\| v-u\|^2_{\mathcal{M}}+\iota_{ \Fix(T)}(v),
$$
where $\iota_{ \Fix(T)}$ is the indicator of $\Fix(T).$
From the first-order optimality conditions (see, e.g., \cite[Theorem 16.47]{Bauschke2017}), we have
\begin{equation*}
0 \in \mathcal{M}(u^* - u) + N_{\Fix(T)}(u^*),
\end{equation*}
which with the definition of the normal cone yields
\begin{equation*}
\langle \mathcal{M}(u - u^*), u^* - v \rangle \geq 0 \quad \forall v \in \Fix(T),
\end{equation*}
or equivalently
\begin{equation}\label{z2}
\langle u - u^*, u^* - v \rangle_{\mathcal{M}} \geq 0, \quad \forall v \in \Fix(T).
\end{equation}
Let \( u^{**} \) be another solution, then it holds
\begin{equation}\label{z3}
\langle u - u^{**}, u^{**} - v \rangle_{\mathcal{M}} \geq 0, \quad \forall v \in \Fix(T).
\end{equation}
Replace \( v \) with \( u^{**} \), \( u^* \) in (\ref{z2}), (\ref{z3}), respectively. This yields
\begin{align*}
\langle u - u^*, u^* - u^{**} \rangle_{\mathcal{M}} &\geq 0,  \\
\langle u - u^{**}, u^{**} - u^* \rangle_{\mathcal{M}} &\geq 0.
\end{align*}
Adding the above two inequalities, we then obtain $\|u^{**} - u^*\|^2_{\mathcal{M}} = 0$. Since
\[
\|u^{**} - u^*\| = \|Tu^{**} - Tu^{*}\| \leq L \|u^{**} - u^*\|_{\mathcal{M}} = 0,
\]
it follows $u^* = u^{**}$, which completes the proof.
\end{proof}

For a point $u$, the $\mathcal{M}$-projection of $u$ onto $\Fix(T)$,  denoted by $P^{\mathcal{M}}_{\Fix(T)}(u)$, is defined by
$$
P^{\mathcal{M}}_{\Fix(T)}(u)=\argmin_{v\in\Fix(T)}\|u-v\|_{\mathcal{M}}.
$$

	
	\begin{lemma}\label{lem24}
		Let $\mathcal{M}:\mathcal{H}\rightarrow \mathcal{H}$ be a linear, bounded, self-adjoint and positive semidefinite operator. Let $T: \mathcal{H} \rightarrow \mathcal{H}$ be  $\mathcal{M}$-nonexpansive  and $\mathcal{M}$-to-norm \(L\)-Lipschitz continuous with a  nonempty closed convex $\Fix({T})$. Then $\forall u \in \mathcal{H}$, the following are equivalent
		\begin{enumerate}
			\item[{\rm(i)}] $u^* = P^{\mathcal{M}}_{\Fix(T)}(u)$;
			\item[{\rm(ii)}] $\langle u^* - u, v - u^* \rangle_{\mathcal{M}} \geq 0, \, \forall v \in \Fix(T)$.
		\end{enumerate}
	\end{lemma}
\begin{proof}
See \cite[Lemma B.1, B.3]{Zhang} and the proof of Lemma \ref{n1}.
\end{proof}
	
	\begin{lemma}
		\label{lem21}
		Let $\mathcal{M} : \mathcal{H}\rightarrow \mathcal{H}$ be a linear, bounded, self-adjoint and positive semidefinite operator and $T: \mathcal{H} \rightarrow \mathcal{H}$ be an $\mathcal{M}$-nonexpansive operator. Let $\{x^k\}$ be a sequence such that  $x^k \rightharpoonup x$ and $\|(I-T)x^k \|_{\mathcal{M}} \rightarrow 0$ as $k\rightarrow \infty $. Then $\|(I-T)x\|_{\mathcal{M}} = 0$.  Furthermore, if $T$ is $\mathcal{M}$-to-norm \(L\)-Lipschitz continuous with  $\Fix({T})\neq \emptyset $, then  $Tx\in \Fix(T)$.
	\end{lemma}
	\begin{proof}
		Suppose $\|(I-T)x\|_{\mathcal{M}} \neq  0$. By Opial property \cite{Opial1967, Naldi2025}, we obtain
		$$
		\begin{aligned}
			\liminf_{k \to \infty} \|x^{k} - x\|_{\mathcal{M}}
			&< \liminf_{k \to \infty} \|x^{k} -  Tx\|_{\mathcal{M}} \\
			&\leq \liminf_{k \to \infty} \left( \|(I-T)x^k \|_{\mathcal{M}} + \|Tx^k - Tx\|_{\mathcal{M}} \right) \\
			&= \liminf_{k \to \infty} \|Tx^k - Tx\|_{\mathcal{M}} \\
			&\leq \liminf_{k \to \infty} \|x^k - x\|_{\mathcal{M}},
		\end{aligned}
		$$
		which yields a contradiction.
		Therefore, $\|(I-T)x\|_{\mathcal{M}} =  0$.  Furthermore, from the $\mathcal{M}$-to-norm \(L\)-Lipschitz continuity  of $T$,  it follows $\|T(Tx) - Tx\| \leq L\|Tx - x\|_{\mathcal{M}}=0$ and hence $Tx\in \Fix (T).$ Thus the proof is completed.
	\end{proof}

\begin{remark}
Note that if  $T=(\mathcal{M} + \mathcal{A})^{-1} \mathcal{M}$, the assumption that $T$ is $\mathcal{M}$-to-norm \(L\)-Lipschitz continuous in Lemma \ref{lem21}  is unnecessary. This follows from $\mathcal{M}Tx = \mathcal{M}x$, which implies:
    \[
    T(Tx) = (\mathcal{M} + \mathcal{A})^{-1} \mathcal{M}Tx = (\mathcal{M} + \mathcal{A})^{-1} \mathcal{M}x = Tx.
    \]
\end{remark}
	
	For the convergence analysis, we  deduce the following lemma.
	\begin{lemma}
		\label{lem27}
		Let $\mathcal{M} : \mathcal{H}\rightarrow \mathcal{H}$ be a linear, bounded, self-adjoint and positive semidefinite operator. Let $T: \mathcal{H} \rightarrow \mathcal{H}$ be  $\mathcal{M}$-nonexpansive  and $\mathcal{M}$-to-norm \(L\)-Lipschitz continuous with  $\Fix({T})\neq \emptyset$. Let $\{x^k\}$ be a  sequence such that
		\begin{enumerate}
			\item[\rm(i)] for all $x^* \in \Fix(T)$, $\lim_{k \to \infty} \| x^k - x^* \|_{\mathcal{M}}$ exists;
			\item[\rm(ii)] $\omega_w(x^k) \subseteq \Fix(T)$;
			\item[\rm(iii)] $\lim_{k \to \infty} \| x^k - {T}x^k \| = 0$.
		\end{enumerate}
		Then $\{x^k\}$ converges weakly to a fixed point of $T$.
	\end{lemma}
	\begin{proof}
		From the $\mathcal{M}$-to-norm \(L\)-Lipschitz continuity  of $T$ and the  condition (i), we can deduce that $\{Tx^k\} $ is bounded. By condition (iii), we obtain that $\{x^k\} $ is also bounded.
		Therefore,
		$\{x^k\} $  has at least one weak cluster point.  By (ii),  we only need to show that the weak cluster point is unique. To this end, let $x$ and $y$ be weak sequential cluster points of $\{x^k\} $ and assume  $\{x^{k_i}\}\rightharpoonup x $ and $\{x^{k_j}\} \rightharpoonup y $. Due to $x,y\in\Fix(T)$, the sequences $\| x^k - x \|_{\mathcal{M}}$ and $\| x^k - y \|_{\mathcal{M}}$ converge. Since
		\begin{equation*}\label{h1}
			2\langle x^k , x - y \rangle_{\mathcal{M}} = \| x^k - y \|_{\mathcal{M}}^2 - \| x^k - x \|_{\mathcal{M}}^2 + \| x \|_{\mathcal{M}}^2 - \| y \|_{\mathcal{M}}^2,
		\end{equation*}
		$\{\langle x^k , x - y \rangle_{\mathcal{M}}\} $ converges as well, say $\langle x^k , x - y \rangle_{\mathcal{M}}\rightarrow l$. Passing to the limit along $\{x^{k_i}\}$ and along $\{x^{k_j}\}$ yields, respectively, $l= \langle x , x - y \rangle_{\mathcal{M}}=\langle y , x - y \rangle_{\mathcal{M}}$. Thus, $\| x - y \|_{\mathcal{M}}=0$. By the $\mathcal{M}$-to-norm \(L\)-Lipschitz continuity of $T$, we have $\|Tx - Ty\|=0$. Specifically, $x =Tx =Ty =y$. We complete the proof.
	\end{proof}
	\begin{remark}\label{rem}
		The assumption (ii) of Lemma \ref{lem27} can be removed when  $T=(\mathcal{M} + \mathcal{A})^{-1}\mathcal{M}$ where $\mathcal{A} : \mathcal{H} \to 2^{\mathcal{H}}$ is a maximal monotone operator with $\zer(\mathcal{A}) \neq \emptyset$, and $\mathcal{M}$ is an admissible preconditioner for $\mathcal{A}$ such that $(\mathcal{M} + \mathcal{A})^{-1}$ is $L$-Lipschitz continuous.
		Indeed,  since $\mathcal{M}(x^k-Tx^k)\in \mathcal{A}(Tx^k)$, by Lemma \ref{lem27}(iii)  and the weak-strong closed property of maximal monotone operators,  we have $\omega_w(Tx^k) \subseteq\Fix(T)$. Using Lemma \ref{lem27}(iii) again, we obtain $\omega_w(x^k) \subseteq\Fix(T)$.
	\end{remark}

	\section{A preconditioned  Halpern iteration with adaptive anchoring parameters}

	In this section, we introduce a preconditioned  Halpern iteration with adaptive anchoring parameters (PHA) and establish the convergence properties of the proposed algorithm. Furthermore, we extend the convergence results to Halpern-type preconditioned proximal point method with adaptive anchoring parameters.
	
	Let $\mathcal{M} : \mathcal{H}\rightarrow \mathcal{H}$ be a linear, bounded, self-adjoint and positive semidefinite operator and $T: \mathcal{H} \rightarrow \mathcal{H}$ be an $\mathcal{M}$-nonexpansive operator. Now  we extend Halpern iteration with adaptive anchoring parameters  in \cite{He} to propose a method to approximate a fixed point of an  $\mathcal{M}$-nonexpansive operator as follows.
	\begin{algorithm} 
		\label{Alg31} 
		\rm{
			\hrule
			\noindent\textbf{\footnotesize{Preconditioned  Halpern iteration with adaptive anchoring parameters (PHA).}}
			\hrule
			\vskip 1mm
			
			\noindent
			1. {\bf Input:}
			initial point $x^0 \in \mathcal{H}$.
			
			\noindent
			2. {\bf for} $k =  1, \ldots$ {\bf do}
			
			\noindent
			3.\quad
			\begin{equation}
				\label{a1}
				x^{k} = \frac{1}{\varphi_k+1}x^{0} + \frac{\varphi_k}{\varphi_k+1}Tx^{k-1},
			\end{equation}
			\quad where
			\begin{equation}
				\label{a2}
				\varphi_k = \frac{2\langle x^{k-1}-Tx^{k-1},x^{0}-x^{k-1}\rangle_{\mathcal{M}}}{\|x^{k-1} - Tx^{k-1}\|^2_{\mathcal{M}}} + 1.
			\end{equation}
			
			\noindent
			4. \textbf{end \ for}
			\vskip 1mm
			\hrule
			
			\hspace*{\fill}
		}
	\end{algorithm}

	Algorithm \ref{Alg31} becomes Halpern-type preconditioned proximal point method with adaptive anchoring parameters when $T=(\mathcal{M} + \mathcal{A})^{-1}\mathcal{M}$   where $\mathcal{A} : \mathcal{H} \to 2^{\mathcal{H}}$ is a maximal monotone operator, and $\mathcal{M}$ is an admissible preconditioner of $\mathcal{A}$.

Restart strategies for Halpern-type iterations have recently received increasing attention. By periodically or adaptively updating the anchoring points and anchoring parameters, such techniques have been shown to improve convergence speed and robustness in applications ranging from image processing to machine learning \cite{Zhang,Park2022,Chen2024}. In particular, restarted Halpern iterations incorporated into PDHG have demonstrated accelerated convergence properties \cite{Park2022,Lu2024}. Motivated by these results, we propose a restarted version of the preconditioned Halpern iteration. We denote by $x^k=\textrm{PHA}(x^{k-1};x^0)$ the update defined in \eqref{a1}-\eqref{a2}, and describe the restarted scheme in  Algorithm \ref{Alg32}.

\begin{algorithm}
    \label{Alg32}
    \rm{
        \hrule
        \noindent\textbf{\footnotesize{Restarted PHA.}}
        \hrule
        \vskip 1mm
        \noindent
        \,\quad\textbf{Input:}
        Initial point \( x^{0,0} \), outer loop counter $n \gets 0$.

        \noindent
        \begin{tabular}{@{}c@{\hspace{1mm}}l@{}}
            \multicolumn{2}{@{}l@{}}{1.\textbf{ repeat}} \\[2pt]
          2.  \hspace{0.5mm} \vline & \quad initialize the inner loop counter $k \gets 1$; \\
          3.  \hspace{0.5mm} \vline & \quad \textbf{repeat} \\
          4.  \hspace{0.5mm} \vline & \quad \, \vline \hspace{3mm} $ x^{n,k} \gets \text{PHA} (x^{n,k-1}; x^{n,0})$; \\
          5.  \hspace{0.5mm} \vline & \quad \textbf{until} restart condition holds; \\
          6.  \hspace{0.5mm} \vline & \quad initialize the initial solution $x^{n+1,0} \gets x^{n,k}$; \\
          7.  \hspace{0.5mm} \vline & \quad $n \gets n+1 $; \\[2pt]
            \multicolumn{2}{@{}l@{}}{8.\textbf{ until} {$x^{n+1,0}$ converges;}}
        \end{tabular}

        \vskip 1mm
        \hrule

        \hspace*{\fill}
    }
\end{algorithm}

	We now  establish some properties of the sequence $\{\varphi_k\}$, whose proof can be given by following the arguments in the proof of {\rm \cite[Lemma 3.1]{He}}.
	
	\begin{lemma} \label{lem26}
		The following hold for all $k \geq 1$,
		\begin{enumerate}
			\item[{\rm(i)}] $\varphi_k \geq k$,
			\item[{\rm(ii)}] $\|x^k - T x^k\|^2_{\mathcal{M}} \leq \frac{2}{\varphi_k} \langle x^k - T x^k, x^0 - x^k \rangle_{\mathcal{M}}.$
		\end{enumerate}
	\end{lemma}
	
	In the following, we present the convergence results of Algorithm \ref{Alg31}.
	\begin{theorem}
		\label{th31}
		Let $\mathcal{M} : \mathcal{H}\rightarrow \mathcal{H}$ be a linear, bounded, self-adjoint and positive semidefinite operator. Let $T: \mathcal{H} \rightarrow \mathcal{H}$ be  $\mathcal{M}$-nonexpansive  and $\mathcal{M}$-to-norm \(L\)-Lipschitz continuous with  a  nonempty closed convex $\Fix({T})$.  Let $\{x^k\}$ be the sequence generated by Algorithm \ref{Alg31} such that $\omega_w(x^k)\subseteq\Fix(T)$. Then $\{x^k\}$ converges strongly to a fixed point of ${T}$.
	\end{theorem}
	
	\begin{proof}
		Firstly, we show that $\{x^{k}\}$ is bounded. Take arbitrarily $p \in \Fix({T})$. From \eqref{a1}, we get
		\begin{equation}
			\label{eq1}
			\begin{aligned}
				\|x^{k} - p\|_{\mathcal{M}}
				&= \left\|\frac{1}{\varphi_k+1}(x^0 - p) + \frac{\varphi_k}{\varphi_k+1}(Tx^{k-1} - p)\right\|_{\mathcal{M}} \\
				&\leq \frac{1}{\varphi_k+1}\|x^0 - p\|_{\mathcal{M}} + \frac{\varphi_k}{\varphi_k+1} \|Tx^{k-1} - p\|_{\mathcal{M}} \\
				&\leq \frac{1}{\varphi_k+1}\|x^0 - p\|_{\mathcal{M}} + \frac{\varphi_k}{\varphi_k+1} \|x^{k-1} - p\|_{\mathcal{M}} \\
				&\leq \max\{\|x^0 - p\|_{\mathcal{M}}, \|x^{k-1} - p\|_{\mathcal{M}}\}.
			\end{aligned}
		\end{equation}
		By induction, we obtain
		\begin{equation}
			\label{eq2}
			\|x^{k} - p\|_{\mathcal{M}} \le\|x^0 - p\|_{\mathcal{M}}.
		\end{equation}
		The  $\mathcal{M}$-to-norm \(L\)-Lipschitz continuity   of $T$ and \eqref{eq2} together give
		\begin{equation}
			\label{eq3}
			\begin{aligned}
				\|x^{k} - p\|
				&\leq \frac{1}{\varphi_k+1}\|x^0 - p\| + \frac{\varphi_k}{\varphi_k+1} \|Tx^{k-1} - p\| \\
				&\leq \frac{1}{\varphi_k+1}\|x^0 - p\| + \frac{\varphi_k}{\varphi_k+1}L\|x^{k-1} - p\|_{\mathcal{M}} \\
				&\leq \max\{\|x^0 - p\|, L\|x^{k-1} - p\|_{\mathcal{M}}\}\\
				&\le\max\{\|x^0 - p\|, L\|x^{0} - p\|_{\mathcal{M}}\},
			\end{aligned}
		\end{equation}
		which  yields that $\{x^{k}\}$ is bounded and hence $\omega_w(x^k)\neq \emptyset$. It follows from the $\mathcal{M}$-to-norm \(L\)-Lipschitz continuity of $T$ that the sequence $\{Tx^k\}$ is bounded too.
		Now, we shall show  $\lim_{k \rightarrow \infty}\|Tx^k - x^k\|_{\mathcal{M}}= 0$.
		By Lemma \ref{lem26}, it holds
		\begin{equation}
			\begin{aligned}
				\label{a13}
				\|x^{k} - Tx^k\|_{\mathcal{M}}\leq\frac{2}{\varphi_k} \|x^{0} - x^k\|_{\mathcal{M}}.
			\end{aligned}
		\end{equation}
		So, the boundedness of $\{x^{k}\}$  ensures that
		\begin{equation}
			\label{a14}
			\lim_{k \rightarrow \infty}\|Tx^k - x^k\|_{\mathcal{M}} = 0 .
		\end{equation}
		
		Next we distinguish two cases to complete the proof.
		
		\vspace{1em}
		\noindent
		Case 1. $\sum_{k=0}^{\infty} \frac{1}{\varphi_k+1}= \infty$.
		\vspace{1em}
		
		In this case, we shall prove that $\lim_{k \rightarrow \infty}\|x^k - q\| = 0$, where
		$
		q = P_{\text{Fix}({T})}^{\mathcal{M}}(x^0).
		$
		Firstly, we show
		\[
		\limsup_{k \to \infty} \langle x^0 - q, x^k - q \rangle_{\mathcal{M}} \leq 0.
		\]
		We take $\{x^{k_j}\} \subset \{x^{k}\}$ such that
		\[
		\limsup_{k \to \infty} \langle x^0 - q, x^k - q \rangle_{\mathcal{M}} = \lim_{j \to \infty} \langle x^0 - q, x^{k_j} - q \rangle_{\mathcal{M}}.
		\]
		Without loss of generality, we  assume \( x^{k_j} \rightharpoonup \bar{x}\) as  $j \to \infty$. From Lemma \ref{lem21} and (\ref{a14}), we  get $T\bar{x}\in \Fix(T)$ and $\|\bar{x}-T\bar{x}\|_{\mathcal{M}}=0$.  By Lemma \ref{lem24} and  Cauchy--Schwarz inequality with respect to the $\mathcal{M}$-scalar product, we have
		\begin{equation}
			\label{b1}
			\begin{aligned}
				\limsup_{k \to \infty} \langle x^0 - q, x^k - q \rangle_{\mathcal{M}}&=\langle x^0 - q, \bar{x} - q \rangle_{\mathcal{M}}\\
				&= \langle x^0 - q, T\bar{x} - q \rangle_{\mathcal{M}}+\langle x^0 - q, \bar{x}-T\bar{x}\rangle_{\mathcal{M}}\\
				&\leq \| x^0 - q\|_{\mathcal{M}}\| \bar{x}-T\bar{x}\|_{\mathcal{M}}=0.
			\end{aligned}
		\end{equation}
		From \eqref{a1}, it follows
		\begin{equation}
			\label{f1}
			\begin{aligned}
				&\|x^{k} - q\|^2_{\mathcal{M}}\\
				=& \left\|\frac{1}{\varphi_k+1}(x^0 - q) + \frac{\varphi_k}{\varphi_k+1}(Tx^{k-1} - q)\right\|^2_{\mathcal{M}}\\
				=& \frac{1}{(\varphi_k+1)^2}\|x^0 - q\|^2_{\mathcal{M}} + \frac{\varphi_k^2}{(\varphi_k+1)^2} \|Tx^{k-1} - q\|^2_{\mathcal{M}} + \frac{2\varphi_k}{(\varphi_k+1)^2} \langle x^0 - q, Tx^{k-1} - q \rangle_{\mathcal{M}}\\
				\leq& \frac{1}{(\varphi_k+1)^2} \|x^0 - q\|^2_{\mathcal{M}} + \frac{\varphi_k}{\varphi_k+1} \|x^{k-1} - q\|^2_{\mathcal{M}} + \frac{2\varphi_k}{(\varphi_k+1)^2} \langle x^0 - q, x^{k-1} - q \rangle_{\mathcal{M}}\\
				&+\frac{2\varphi_k}{(\varphi_k+1)^2} \|x^0 - q\|_{\mathcal{M}}\|Tx^{k-1} - x^{k-1}\|_{\mathcal{M}}.
			\end{aligned}
		\end{equation}
		Let \( a_{k+1} = \|x^{k} - q\|^2_{\mathcal{M}}\), \( \gamma_k = \frac{1}{\varphi_k+1} \) and
		\[
		\delta_k = \frac{1}{\varphi_k+1} \|x^0 - q\|^2_{\mathcal{M}} + \frac{2\varphi_k}{\varphi_k+1} ( \langle x^0 - q, x^{k-1} - q \rangle_{\mathcal{M}} + \|x^0 - q\|_{\mathcal{M}} \|Tx^{k-1} - x^{k-1}\|_{\mathcal{M}} ).
		\]
		Then we can rewrite (\ref{f1})  as follows
		\begin{equation*}
			a_{k+1} \leq (1 - \gamma_k)a_k + \gamma_k \delta_k, \quad k \geq 0.
		\end{equation*}
		The facts $\varphi_k \to \infty$ and (\ref{a14}) together with (\ref{b1}) imply $\limsup_{k \to \infty} \delta_k \leq 0$.
		We  have $\sum_{k=0}^{\infty} \gamma_k = \infty$ from the assumption of Case 1.
		By \cite[Lemma 2.5]{Xu}, we conclude
		\begin{equation}
			\label{f2}
			\lim_{k \to \infty} \|x^{k} - q\|_{\mathcal{M}} = 0.
		\end{equation}
		By \eqref{eq3}, we have
		\begin{equation*}
			\begin{aligned}
				\|x^{k} - q\|
				&\leq \frac{1}{\varphi_k+1}\|x^0 - q\| + \frac{\varphi_k}{\varphi_k+1}L\|x^{k-1} - q\|_{\mathcal{M}}.
			\end{aligned}
		\end{equation*}
		It easily follows that $\lim_{k \rightarrow \infty}\|x^{k} - q\|= 0$ from (\ref{f2}) and $\frac{1}{\varphi_k+1} \to 0$ as $k \to \infty$ . This finishes the proof of Case 1.
		
		\vspace{1em}
		\noindent
		Case 2. $\sum_{k=0}^{\infty} \frac{1}{\varphi_k+1} < \infty $.
		\vspace{1em}
		
		Set \( \gamma_k = \frac{1}{\varphi_k+1} \). By the assumption, it holds \( \sum_{k=0}^{\infty} \gamma_k < \infty \).
		Now we  show $ \lim_{k \to \infty} \|x^{k} -Tx^{k}\| =0 $. By (\ref{a1})
		\begin{equation}
			\begin{aligned}
				\label{c10}
				\|x^{k+1} -Tx^{k+1}\|
				=& \|Tx^{k} -Tx^{k+1}+\gamma_{k+1}(x^0 - Tx^k)\| \\
				\leq & L\|x^k -x^{k+1}\|_{\mathcal{M}}+\gamma_{k+1}\|x^0 - Tx^k\|,
			\end{aligned}
		\end{equation}
		where the inequality comes from the $\mathcal{M}$-to-norm \(L\)-Lipschitz continuity  of $T$. Next we verify $ \lim_{k \to \infty} \|x^{k} - x^{k-1}\|_{\mathcal{M}} =0 $.
		Using (\ref{a1}) again, we get
		\begin{equation}
			\begin{aligned}
				\label{b3}
				\|x^k - x^{k-1}\|_{\mathcal{M}}
				&= \left\| \frac{1}{\varphi_k+1} \left(x^0 - x^{k-1}\right) + \frac{\varphi_k}{\varphi_k+1} \left( Tx^{k-1} - x^{k-1} \right) \right\|_{\mathcal{M}} \\
				&\leq \gamma_k \|x^0 - x^{k-1}\|_{\mathcal{M}} + \|x^{k-1} - Tx^{k-1}\|_{\mathcal{M}}.
			\end{aligned}
		\end{equation}
		Rewrite the first inequality in (\ref{a13}) in terms of $\gamma_k$ as
		\begin{equation*}
			\begin{aligned}
				\|x^k - Tx^k\|_{\mathcal{M}} \leq \frac{2\gamma_k}{1 - \gamma_k} \|x^0 - x^k\|_{\mathcal{M}} \leq 4\gamma_k \|x^0 - x^k\|_{\mathcal{M}},
			\end{aligned}
		\end{equation*}
		due to $\gamma_k \leq \frac{1}{k+1} \leq \frac{1}{2}$ for all $k \geq 1$. From the  boundedness of $\{x^k\}$,  it turns out that
		\begin{equation}
			\begin{aligned}
				\label{b2}
				\sum_{k=0}^{\infty} \|x^k - Tx^k\|_{\mathcal{M}} < \infty.
			\end{aligned}
		\end{equation}
		Putting (\ref{b3}) and (\ref{b2})   together leads to
		\begin{equation}
			\begin{aligned}
				\label{b5}
				\sum_{k=0}^{\infty} \|x^k - x^{k-1}\|_{\mathcal{M}} < \infty,
			\end{aligned}
		\end{equation}
		which yields
		\begin{equation}
			\begin{aligned}
				\label{b4}
				\lim_{k \rightarrow \infty}\|x^k - x^{k-1}\|_{\mathcal{M}}= 0.
			\end{aligned}
		\end{equation}
		From (\ref{c10}), (\ref{b4}), the boundedness of $\{Tx^k\}$ and the fact that $\gamma_k \to 0$  as $k \to \infty$, we derive that
		\begin{equation}\label{h3}
			\lim_{k \to \infty}\| x^k-T x^k \|= 0.
		\end{equation}
		By \eqref{eq1}, it follows
		\begin{equation*}
			\label{f3}
			\begin{aligned}
				\|x^k - p\|_{\mathcal{M}}
				&\leq  \gamma_k \|x^0 - p\|_{\mathcal{M}} + (1-\gamma_k) \|Tx^{k-1} - p\|_{\mathcal{M}} \\
				&\leq \|x^{k-1} - p\|_{\mathcal{M}} + \gamma_k \|x^0 - p\|_{\mathcal{M}},
			\end{aligned}
		\end{equation*}
		which with \cite[Lemma 2 in Chapter 2]{Polyak} yields that
		\( \lim_{k \to \infty} \|x^k - p\|_{\mathcal{M}} \) exists.
		Combining (\ref{h3}) and Lemma \ref{lem27}, we get that $\{x^k\}$ converges weakly to a fixed point of $T$.
		
		Now we  show that $\{x^k\}$ is a Cauchy sequence. Due to the boundedness of $\{Tx^k\}$, there exists $M>0$ such that $\|x^0 - Tx^{k}\|\leq M$ for $k\ge0$. By (\ref{a1})
		\begin{equation*}
			\begin{aligned}
				\label{c3}
				\|x^k - x^{k-1}\|
				=&\|\gamma_k x^0 + (1-\gamma_k)Tx^{k-1}-\gamma_{k-1} x^0 - (1-\gamma_{k-1})Tx^{k-2}\| \\
				\leq &(1-\gamma_k)\|Tx^{k-1} - Tx^{k-2}\|+|\gamma_k-\gamma_{k-1}|\|x^0 - Tx^{k-2}\|\\
				\leq &L(1-\gamma_k)\|x^{k-1} - x^{k-2}\|_{\mathcal{M}}+M(\gamma_k+\gamma_{k-1}).
			\end{aligned}
		\end{equation*}
		Combining (\ref{b5}) and $\sum_{k=1}^{\infty} \gamma_k < \infty$ yields that $\sum_{k=1}^{\infty}\|x^k - x^{k-1}\| < \infty$ and  $\{x^k\}$ is a Cauchy sequence. Therefore, $\{x^k\}$ converges strongly to a fixed point of $T$. The proof is  complete.
	\end{proof}

 As the Halpern
iteration with adaptive anchoring parameters in \cite{He} speeds up   the standard PPP algorithm, Algorithm \ref{Alg31} can be used to  accelerate   the degenerate PPP algorithm and the splitting methods it encompasses.
	
	By (\ref{h3}), Theorem \ref{th31},  Lemma \ref{lem211}  and Remark \ref{rem}, we obtain the following convergence result for Halpern-type preconditioned proximal point method with adaptive anchoring parameters.

	\begin{corollary}
		\label{c1}
		Let $\mathcal{A} : \mathcal{H} \to 2^{\mathcal{H}}$ be a maximal monotone operator with $\zer(\mathcal{A}) \neq \emptyset$, and $\mathcal{M}$ be an admissible preconditioner for $\mathcal{A}$ such that $(\mathcal{M} + \mathcal{A})^{-1}$ is Lipschitz continuous. Let $\{x^k\}$ be the sequence generated by Algorithm \ref{Alg31} with $T=(\mathcal{M} + \mathcal{A})^{-1}\mathcal{M}$. Then $\{x^k\}$ converges strongly to a point of $\zer(\mathcal{A})$.
	\end{corollary}
	
	Below we present the convergence rate of Algorithm \ref{Alg31}.
	
	\begin{theorem}
		\label{l1}
		Let $\mathcal{M} : \mathcal{H}\rightarrow \mathcal{H}$ be a linear, bounded, self-adjoint and positive semidefinite operator. Let $T: \mathcal{H} \rightarrow \mathcal{H}$ be $\mathcal{M}$-nonexpansive  and $\mathcal{M}$-to-norm \(L\)-Lipschitz continuous with  $\Fix({T})\neq \emptyset$. Then for any sequence $\{x^k\}$ generated by Algorithm \ref{Alg31}, we have
		\begin{enumerate}
			\item[{\rm(i)}]
			$
			\|x^k - Tx^k\|_{\mathcal{M}} \leq \frac{2}{\varphi_k + 1} \|x^0 - x^*\|_{\mathcal{M}},\,\, k \geq 1,
			$
			\item[{\rm(ii)}]
			$
			\|x^k - Tx^k\| \leq \frac{M(L+1)}{\varphi_k+1}+\frac{2L}{\varphi_{k-1} + 1}\|x^0 - x^*\|_{\mathcal{M}},\,\, k \geq 2,
			$
		\end{enumerate}
		where $M= \sup_{k\ge1}\max\{\|x^0 - x^k\|_{\mathcal{M}}, \|x^0 - Tx^k\|\}$, $\varphi_k$ is defined in (\ref{a2}) and $x^*$ is an arbitrary fixed point of $T$.
	\end{theorem}
	\begin{proof}
		(i) can be proven by employing arguments which are similar to those used in the proof of Theorem 3.2 in \cite{He}.
		
		(ii) By (\ref{b3}) and  (i), we have
		\begin{equation}
			\begin{aligned}
				\label{b333}
				\|x^k - x^{k-1}\|_{\mathcal{M}}
				&\leq \frac{1}{\varphi_k+1}\|x^0 - x^{k-1}\|_{\mathcal{M}} + \|x^{k-1} - Tx^{k-1}\|_{\mathcal{M}}\\
				&\leq \frac{1}{\varphi_k+1}M+\frac{2}{\varphi_{k-1} + 1}\|x^0 - x^*\|_{\mathcal{M}}.
			\end{aligned}
		\end{equation}
		From (\ref{c10}) and (\ref{b333}), we derive
		\begin{equation*}
			\begin{aligned}
				\label{d1}
				\|x^{k} -Tx^{k}\|
				\leq & L\|x^{k-1} -x^{k}\|_{\mathcal{M}}+\frac{1}{\varphi_k+1}\|x^0 - Tx^{k-1}\|\\
				\leq & \frac{M(L+1)}{\varphi_k+1}+\frac{2L}{\varphi_{k-1} + 1}\|x^0 - x^*\|_{\mathcal{M}},
			\end{aligned}
		\end{equation*}
		which completes the proof.
	\end{proof}
	
	Similar to Corollary \ref{c1}, we  get the convergence rate of Halpern-type preconditioned proximal point method with adaptive anchoring parameters by using Theorem \ref{l1} and  Lemma \ref{lem211}.
	\begin{corollary}
		\label{c22}
		Let $\mathcal{A} : \mathcal{H} \to 2^{\mathcal{H}}$ be a maximal monotone operator with $\zer(\mathcal{A}) \neq \emptyset$, and $\mathcal{M}$ be an admissible preconditioner for $\mathcal{A}$ with a decomposition $\mathcal{M}=\mathcal{C}\mathcal{C}^*$ such that $(\mathcal{M} + \mathcal{A})^{-1}$ is $L$-Lipschitz continuous. Let $\{x^k\}$ be the sequence generated by Algorithm \ref{Alg31} with $T=(\mathcal{M} + \mathcal{A})^{-1}\mathcal{M}$, then we have
		\begin{enumerate}
			\item[{\rm(i)}]
			$
			\|x^k - Tx^k\|_{\mathcal{M}} \leq \frac{2}{\varphi_k + 1} \|x^0 - x^*\|_{\mathcal{M}},\,\, k \geq 1,
			$
			\item[{\rm(ii)}]
			$
			\|x^k - Tx^k\| \leq \frac{M(L'+1)}{\varphi_k+1}+\frac{2L'}{\varphi_{k-1} + 1}\|x^0 - x^*\|_{\mathcal{M}},\,\, k \geq 2,
			$
		\end{enumerate}
		where $M= \sup_{k\ge1}\max \{\|x^0 - x^k\|_{\mathcal{M}}, \|x^0 - Tx^k\|\}$, $L' =L \|\mathcal{C}\|$, $\varphi_k$ is defined in (\ref{a2}) and $x^*$ is an arbitrary  point of $\zer({\mathcal{A}})$.
	\end{corollary}

By Corollary \ref{c22}(ii), Halpern-type preconditioned proximal point method with adaptive anchoring parameters achieves a convergence rate of  $\mathcal{O}(1/\varphi_{k})$. By $\varphi_{k}\ge k$ from Lemma  \ref{lem26}(i),  it is a further improvement of the $\mathcal{O}(1/k)$ rate of the HPPP algorithm in \cite{Zhang}. In practice, $\varphi_{k}$ does not grow linearly most of the time (see, e.g.,  numerical experiments in Section \ref{sect5}).
	
	\section{An acceleration of  Chambolle--Pock algorithm}
	In this section, we introduce an acceleration of  Chambolle--Pock algorithm by combining the PHA iteration  and Chambolle--Pock algorithm. Then we establish the strong convergence and convergence rate of the proposed algorithm.
	
	The Chambolle--Pock (CP) algorithm \cite{CP} which is also known as the Primal-Dual Hybrid Gradient (PDHG) method, is an effective method for solving problem \eqref{tv1} addressed as
	\begin{equation}
		\label{CP}
		\left\{
		\begin{aligned}
			{u}_{k+1} &= \prox_{\tau f}\left({u}_k - \tau \mathbf{K}^* v_k\right), \\
			{v}_{k+1} &= \prox_{\sigma g^*}\left(v_k + \sigma \mathbf{K}(2u_{k+1} - u_k)\right),
		\end{aligned}
		\right.
	\end{equation}
	where $\tau, \sigma > 0 $. In the original paper the convergence of  \eqref{CP} was ensured when \( \tau \sigma\| {\bf K} \|^{2} < 1\).
	Since its inception, the CP algorithm \eqref{CP} has received much attention due to its success in solving large scale structured convex optimization  problems.
	
	Let
	\begin{equation}\label{s1}
		\mathcal{A} = \begin{pmatrix}
			\partial f & \mathbf{K}^* \\
			-\mathbf{K} & \partial g^{*}
		\end{pmatrix},\quad
		\mathcal{M} = \begin{pmatrix}
			\frac{1}{\tau}I & -{\bf K}^{*} \\
			-{\bf K} & \frac{1}{\sigma}I
		\end{pmatrix}
	\end{equation} with
	$\tau, \sigma > 0 $. From \cite{Bredies}, $\mathcal{A}$ is a maximal monotone operator and $ \mathcal{M}$ is a  linear, bounded, self-adjoint semidefinite  operator and therefore a preconditioner when $\tau \sigma \| \mathbf{K} \|^2 \le 1$. Recall that
	\begin{equation*}
		(\mathcal{M}+\mathcal{A})^{-1}:(u,v)\mapsto (\prox_{\tau f}(\tau u), \prox_{\sigma g^*} (2\sigma\mathbf{K} \prox_{\tau f}(\tau u)+\sigma v)),
	\end{equation*}
	by using the nonexpansiveness of the proximal operator, we get that $(\mathcal{M}+\mathcal{A})^{-1}$ is $L$-Lipschitz with $L=\max\left\{ \sqrt{2}\sigma,\ \tau\sqrt{8\sigma^2 \|\mathbf{K}\|^2 + 1} \right\}$.
	
	Following \cite{He2012},  Bredies et al. \cite{Bredies} recently reformulated  problem \eqref{tv1} into the following fixed point form
	\begin{equation}
		\label{ad1}
		\hbox{find a point} \,\,x\in U \times V\,\, \hbox{such that}\,\, x = (\mathcal{M} + \mathcal{A})^{-1} \mathcal{M}x
	\end{equation}
	where $\mathcal{A}$ and $\mathcal{M}$ are given as above. They wrote the CP algorithm \eqref{CP}  as the degenerate PPP method \eqref{4a} and then established the weak convergence of the CP algorithm \eqref{CP} for the degenerate case $\tau \sigma \| \mathbf{K} \|^2 = 1$.
	
	Based on the equivalence of \eqref{tv1} and \eqref{ad1}, we apply the PHA method to solve  problem \eqref{ad1} and propose the following accelerated Chambolle--Pock  (aCP)  algorithm.
	\begin{algorithm}
		\label{Alg81}
		\rm{
			\hrule\hrule
			\vskip 1mm
			\noindent\textbf{An accelerated Chambolle--Pock (aCP) algorithm }
			\vskip 1mm
			\hrule\hrule

			\vskip 1mm
			
			\noindent
			1. {\bf Input:}
			initial point \(  (u^0, v^0)\in U \times V \) and $\tau, \sigma > 0 $ such that $\tau \sigma \| \mathbf{K} \|^2 \leq 1$.
			
			\noindent
			2. {\bf for} $k = 0, 1, \ldots$ {\bf do}
			
			\noindent
			3.\begin{equation}
				\quad \quad \quad\begin{cases}
					\label{88}
					p^{k} =  \prox_{\tau f}\left( u^{k} - \tau \mathbf{K}^* v^{k}\right),\\
					
					q^{k} = \prox_{\sigma g^*} \left( \sigma \mathbf{K}(2p^k -  u^{k}) + v^{k} \right),
				\end{cases}
			\end{equation}
			\noindent
			4.\begin{equation*}
				\begin{cases}
					\label{89}
					u^{k+1} = \frac{1}{\varphi_{k}+1} u^{0} + \frac{\varphi_{k}}{\varphi_{k}+1} p^{k},\\
					
					v^{k+1} = \frac{1}{\varphi_{k}+1} v^{0} + \frac{\varphi_{k}}{\varphi_{k}+1} q^{k},
				\end{cases}
			\end{equation*}
			\quad where
			\begin{equation}
				\label{98}
				\varphi_{k} = \frac{
					2 \langle x^{k} - y^{k}, x^0 - x^{k} \rangle_{\mathcal{M}}}{\|x^{k} - y^{k}\|_{\mathcal{M}}^2} + 1,
			\end{equation}
\quad and $x^k = (u^k, v^k)$ { and } $y^k = (p^k , q^k).$
			\vskip 1mm
			\noindent
			5. \textbf{end \ for}
			\vskip 1mm
}
		\hrule\hrule
		\hspace*{\fill}
	\end{algorithm}

We now make the  blanket assumptions on problem (\ref{tv1}).

\begin{assumption}
\label{ass4.1}
 Assume that the set of solutions of (\ref{tv1}) is nonempty and  $0 \in \mathrm{sri}(\dom g - \mathbf{K}(\dom f))$.
\end{assumption}

Under Assumption \ref{ass4.1}, the set of solutions of  problem (\ref{tv2}) is nonempty (see, e.g., \cite[Theorem 16.47]{Bauschke2017}.
It follows from  \cite[Corollaries 28.2.2 and 28.3.1]{Rockafellar1970} that \(u^* \in U\) is a solution of (\ref{tv1}) if and only if there exists \(v^* \in V\) such that \((u^*, v^*)\) is a saddle-point of \(\mathcal{L}(u, v)\), i.e., \(\mathcal{L}(u^*, v) \leq \mathcal{L}(u^*, v^*) \leq \mathcal{L}(u, v^*)\) for all  $(u, v) \in U \times V $.  We denote the set of solutions of the saddle-point problem (\ref{saddle}) by
	\begin{equation}\label{331}
		\Omega=\{(u^*,v^*) \in U \times V : 0 \in \partial f(u^*)+\mathbf{K}^*v^*,\, 0 \in \partial g^*(v^*)- \mathbf{K}u^*\}.
	\end{equation}
	Then $\Omega=\Fix((\mathcal{M} + \mathcal{A})^{-1} \mathcal{M})$.
	
	By Corollaries \ref{c1} and \ref{c22}, we obtain the strong convergence and convergence rate of Algorithm \ref{Alg81}.
	\begin{corollary}
		Let $\tau$ and $ \sigma$ be positive such that $\tau \sigma \| \mathbf{K} \|^2 \leq 1$.
		Let $\{(u^k, v^k)\}$ be the sequence generated by Algorithm \ref{Alg81}. Then $\{(u^k, v^k)\}$  strongly converges to a point  of $\Omega.$
	\end{corollary}

\begin{corollary}
		\label{c33}
		Let   $\mathcal{M}=\mathcal{C}\mathcal{C}^*$
		be a decomposition of $\mathcal{M}$ in \eqref{s1} and  $\tau$, $ \sigma$ be positive such that $\tau \sigma \| \mathbf{K} \|^2 \leq 1$. Let $\{x^k\}$ and $\{y^k\}$ be the sequences generated by Algorithm \ref{Alg81}, then we have
		\begin{enumerate}
			\item[{\rm(i)}]
			$
			\|x^k - y^k\|_{\mathcal{M}} \leq \frac{2}{\varphi_{k-1} + 1} \|x^0 - x^*\|_{\mathcal{M}},\,\, k \geq 2,
			$
			\item[{\rm(ii)}]
			$
			\|x^k - y^k\| \leq \frac{M(L+1)}{\varphi_{k-1}+1}+\frac{2L}{\varphi_{k-2} + 1}\|x^0 - x^*\|_{\mathcal{M}},\,\, k \geq 3,
			$
		\end{enumerate}
		where $M= \sup_{k\ge1}\max \{\|x^0 - x^k\|_{\mathcal{M}}, \|x^0 - y^k\|\}$, $L=\max\left\{ \sqrt{2}\sigma,\ \tau\sqrt{8\sigma^2 \|\mathbf{K}\|^2 + 1} \right\}\|\mathcal{C}\|$, $\varphi_k$ is defined in (\ref{98}) and $x^*$ is a point of $\Omega.$
	\end{corollary}

\begin{remark}
The above convergence results are  applicable to the general inclusion problem \eqref{tv2} with  maximal monotone operators \( A_1 \) and \( A_2 \) instead of  the specific case of \( \partial f \) and \( \partial g \).
\end{remark}

	To get the convergence rate of Algorithm \ref{Alg81}, we begin by considering the residual mapping associated with $\Omega$ as follows
	\begin{equation}
		\label{kkt}
		R(x) := \begin{pmatrix}
			u - \prox_f (u -  \mathbf{K}^* v) \\v - \operatorname{prox}_{g^*} (v +  \mathbf{K} u)
		\end{pmatrix},\quad \forall x=(u, v) \in U \times V.
	\end{equation}
	It is obvious that  $x^*=(u^*,v^*) \in  \Omega $ if and only if $R(x^*) = 0$.
	On the other hand, $(u^*,v^*) \in  \Omega $ is also equivalent to
	\begin{equation*}
		\begin{aligned}
			P(u)= f(u)-f(u^*) + \left\langle \mathbf{K}^*v^*,u-u^*\right\rangle \geq 0, \quad \forall u \in U, \\
			D(v)= g^*(v)-g^*(v^*) - \left\langle \mathbf{K}u^*,v-v^*\right\rangle \geq 0, \quad \forall v \in V.
		\end{aligned}
	\end{equation*}
	Therefore, we can define the following primal-dual gap function, as introduced in \cite{Chang2021},
	\begin{equation}\label{pd}
		G(u,v):=P(u)+D(v), \quad \forall (u, v) \in U \times V.
	\end{equation}

	The next lemma provides the upper bound for the primal-dual gap function $G$ at $(p^k,q^k)$.
	
	\begin{lemma}
		\label{lem41}
		Let \(\{(p^k, q^k)\}\) be the sequence generated by Algorithm \ref{Alg81} and  \((u^*, v^*)\) be a point of $\Omega$. Then for all \(k \geq 0\), we have
		\begin{equation}
			\label{pdgf}
			G(p^k, q^k) \leq \left\langle  p^k-u^*,  \frac{u^{k} - p^{k}}{\tau} + \mathbf{K}^*(q^k-v^k) \right\rangle + \left\langle  q^k -v^* , \frac{v^{k} - q^{k}}{\sigma}+ \mathbf{K}(p^k -  u^{k}) \right\rangle.
		\end{equation}
	\end{lemma}
	\begin{proof}
		From \eqref{d10} and  \eqref{88}, it follows
		\[
		\frac{u^{k} - p^{k}}{\tau}-\mathbf{K}^*v^k  \in \partial f(p^k), \quad \frac{v^{k} - q^{k}}{\sigma}+ \mathbf{K}(2p^k -  u^{k}) \in \partial g^*(q^k).
		\]
		By {\rm \cite[Proposition 16.10]{Bauschke2017}}, we can get
		\[
		\begin{cases}
			f^*(\frac{u^{k} - p^{k}}{\tau}-\mathbf{K}^*v^k) = \langle p^k, \frac{u^{k} - p^{k}}{\tau}-\mathbf{K}^*v^k \rangle - f(p^k), \\
			g(\frac{v^{k} - q^{k}}{\sigma}+\mathbf{K}(2p^k -  u^{k})) = \langle q^k,\frac{v^{k} - q^{k}}{\sigma}+ \mathbf{K}(2p^k -  u^{k}) \rangle - g^*(q^k).
		\end{cases}
		\]
		Summing them up, we  obtain
		\begin{equation}\label{3.22}
			\begin{aligned}
				&f^*\left(\frac{u^{k} - p^{k}}{\tau}-\mathbf{K}^*v^k\right) + g\left(\frac{v^{k} - q^{k}}{\sigma}+\mathbf{K}(2p^k - u^{k})\right) \\
				= &-f(p^k) - g^*(q^k) + \left\langle p^k, \frac{u^{k} - p^{k}}{\tau}-\mathbf{K}^*v^k \right\rangle  + \left\langle q^k, \frac{v^{k} - q^{k}}{\sigma} + \mathbf{K}(2p^k - u^{k}) \right\rangle.
			\end{aligned}
		\end{equation}
		On the other hand, from \eqref{331}, we have
		\begin{equation}\label{3.25}
			u^* \in \partial f^*(-\mathbf{K}^*v^*), \quad  v^* \in \partial g(\mathbf{K}u^*).
		\end{equation}
		Therefore, it follows from the subdifferential inequality, that
		\[
		\begin{cases}
			f^*\left(\frac{u^{k} - p^{k}}{\tau}-\mathbf{K}^*v^k\right) - f^*(-\mathbf{K}^*v^*) \geq \left\langle u^*, \frac{u^{k} - p^{k}}{\tau}-\mathbf{K}^*v^k + \mathbf{K}^*v^* \right\rangle \\
			g\left(\frac{v^{k} - q^{k}}{\sigma}+\mathbf{K}(2p^k - u^{k})\right) - g(\mathbf{K}u^*) \geq \left\langle v^*, \frac{v^{k} - q^{k}}{\sigma}+\mathbf{K}(2p^k - u^{k}) - \mathbf{K}u^* \right\rangle.
		\end{cases}
		\]
		Summing them up, we get
		\begin{equation}\label{3.23}
			\begin{aligned}
				&f^*\left(\frac{u^{k} - p^{k}}{\tau}-\mathbf{K}^*v^k\right) + g\left(\frac{v^{k} - q^{k}}{\sigma}+\mathbf{K}(2p^k - u^{k})\right) \\
				\geq& f^*(-\mathbf{K}^*v^*) + g(\mathbf{K}u^*) + \left\langle u^*, \frac{u^{k} - p^{k}}{\tau}-\mathbf{K}^*v^k + \mathbf{K}^*v^* \right\rangle \\
				&+ \left\langle v^*, \frac{v^{k} - q^{k}}{\sigma}+\mathbf{K}(2p^k - u^{k}) - \mathbf{K}u^* \right\rangle \\
				=&f^*(-\mathbf{K}^*v^*) + g(\mathbf{K}u^*) + \left\langle u^*, \frac{u^{k} - p^{k}}{\tau}-\mathbf{K}^*v^k \right\rangle \\
				&+ \left\langle v^*, \frac{v^{k} - q^{k}}{\sigma}+\mathbf{K}(2p^k - u^{k}) \right\rangle.
			\end{aligned}
		\end{equation}
		Combining \eqref{3.22} and \eqref{3.23}, we deduce
		\begin{equation*}\label{3.24}
			\begin{aligned}
				&f(p^k) + g^*(q^k) + f^*(-\mathbf{K}^*v^*) + g(\mathbf{K}u^*) \\
				\leq & \left\langle p^k-u^*, \frac{u^{k} - p^{k}}{\tau}-\mathbf{K}^*v^k \right\rangle
				+ \langle q^k-v^*,\frac{v^{k} - q^{k}}{\sigma}+ \mathbf{K}(2p^k -  u^{k}) \rangle.
			\end{aligned}
		\end{equation*}
Since  \( f^*(-\mathbf{K}^*v^*) + g(\mathbf{K}u^*) = -f(u^*) - g^*(v^*) \)   from (\ref{3.25}) and {\rm \cite[Proposition 16.10]{Bauschke2017}},
 we derive
\begin{equation}\label{3.2444}
			\begin{aligned}
				&f(p^k) - f(u^*) + g^*(q^k) - g^*(v^*) \\
				\leq & \left\langle p^k-u^*, \frac{u^{k} - p^{k}}{\tau}-\mathbf{K}^*v^k \right\rangle
				+ \langle q^k-v^*,\frac{v^{k} - q^{k}}{\sigma}+ \mathbf{K}(2p^k -  u^{k}) \rangle.
			\end{aligned}
		\end{equation}
It follows from \eqref{pd} that 		
\begin{equation}
\label{3.2444a}
			\begin{aligned}
 f(p^k) - f(u^*) + g^*(q^k) - g^*(v^*) = G(p^k, q^k) - \langle v^*, \mathbf{K}p^k \rangle + \langle u^*, \mathbf{K}^*q^k \rangle.
\end{aligned}
		\end{equation}
Combining \eqref{3.2444} and \eqref{3.2444a} and using $0 = \langle q^k, -\mathbf{K}p^k \rangle + \langle p^k, \mathbf{K}^*q^k \rangle$, we can get \eqref{pdgf}. This completes the proof.
	\end{proof}
	Now, we  present the complexity result of Algorithm \ref{Alg81}.
	\begin{theorem}
		 Let $\mathcal{M}$ and $\mathcal{A}$ be given as in \eqref{s1} and $T=(\mathcal{M} + \mathcal{A})^{-1}\mathcal{M}$. Let $\mathcal{M}=\mathcal{C}\mathcal{C}^*$
		be a decomposition of $\mathcal{M}$  and  $\tau$, $ \sigma$ be positive such that $\tau \sigma \| \mathbf{K} \|^2 \leq 1$.
		Let $ \{(p^k, q^k)\}$  be the sequence generated by Algorithm \ref{Alg81} and  $x^*=(u^*,v^*)$ be a point of  $\Omega$. Then we have the following bounds
		\begin{enumerate}
			\item[{\rm(i)}]
			\begin{equation*}
				\label{92}
				\|R(y^k)\| \leq \frac{M_1 \rho(L+1)}{\varphi_{k-1}+1} + \frac{2\rho L}{\varphi_{k-2} + 1}\|x^0 - x^*\|_{\mathcal{M}}, \quad k \geq 2.
			\end{equation*}
			\item[{\rm(ii)}]
			\begin{equation*}
				G(p^k, q^k) \leq \frac{M_1 M_2 \rho (L+1) }{\varphi_{k-1}+1} + \frac{2 M_2 \rho L}{\varphi_{k-2} + 1}\|x^0 - x^*\|_{\mathcal{M}}, \quad k \geq 2.
			\end{equation*}
		\end{enumerate}
		where $\rho = \sqrt{ 2(\|\mathbf{K}\|^2+\max\left\{ \frac{1}{\tau^2},\,  \frac{1}{\sigma^2} \right\})}$, $L = \max\left\{ \sqrt{2}\sigma,\ \tau\sqrt{8\sigma^2 \|\mathbf{K}\|^2 + 1} \right\} \|\mathcal{C}\|$,  $M_1 = \sup_{k \geq 1}\max \left\{ \|x^0 - x^k\|_{\mathcal{M}}, \|x^0 - Tx^k\| \right\}$,
 $M_2 = \sup_{k \geq 0} \left\{ \|p^k - u^*\| + \|q^k - v^*\| \right\}$ and $\varphi_k$ is defined in \eqref{98}.
		
	\end{theorem}
	\begin{proof}
		From \eqref{d10}  and \eqref{88}, we can derive
		\begin{equation}
			\label{93}
			\begin{cases}
				p^k = \prox_f \left( p^k -  (\frac{p^k - u^k}{\tau}  + \mathbf{K}^* v^k ) \right), \\
				q^k = \prox_{g^*} \left( q^k - ( \frac{q^k - v^k}{\sigma}  - \mathbf{K} (2p^k  - u^k)) \right).
			\end{cases}
		\end{equation}
		By \eqref{93} and the nonexpansiveness of the  proximal operator, we have
		\begin{equation}
			\label{94}
			\begin{aligned}
				\| p^k - \prox_f (p^k - \mathbf{K}^* q^k) \|
				&= \left\| \prox_f \left( p^k - ( \frac{p^k - u^k}{\tau}  + \mathbf{K}^* v^k ) \right) - \prox_f (p^k - \mathbf{K}^* q^k) \right\| \\
				&\leq \left\| -\frac{1}{\tau} (p^k - u^k)-\mathbf{K}^* ( v^k-q^k) \right\| \\
				&\leq   \frac{1}{\tau} \| p^k - u^k \|+ \|\mathbf{K}\| \| v^k - q^k \|.
			\end{aligned}
		\end{equation}
		Similarly, we get
		\begin{equation}
			\label{95}
			\begin{aligned}
				&\| q^k - \prox_{g^*} (q^k + \mathbf{K} p^k) \|\\
				=& \left\| \prox_{g^*} \left( q^k - \left( \frac{q^k - v^k}{\sigma}  - \mathbf{K} (2p^k  - u^k) \right) \right) - \prox_{g^*}(q^k + \mathbf{K} p^k) \right\| \\
				\leq& \left\| -\frac{1}{\sigma} (q^k - v^k) + \mathbf{K} ( p^k - u^k) \right\| \\
				\leq& \frac{1}{\sigma} \| q^k - v^k \|+ \|\mathbf{K}\| \| p^k - u^k \|.
			\end{aligned}
		\end{equation}
		By \eqref{kkt}, \eqref{94}, \eqref{95} and arithmetic-mean inequality,  we can obtain
		\begin{equation}
			\label{96}
			\begin{aligned}
				&\| R(y^k) \|^2\\
				\leq& 2 \left( \frac{1}{\tau^2} \| p^k - u^k \|^2 + \|\mathbf{K}\|^2 \| v^k - q^k \|^2  \right) + 2 \left( \frac{1}{\sigma^2} \| q^k - v^k \|^2 + \|\mathbf{K}\|^2 \| p^k - u^k \|^2 \right) \\
				=& 2(\|\mathbf{K}\|^2 + \frac{1}{\sigma^2}) \| v^k - q^k \|^2 + 2(\|\mathbf{K}\|^2 + \frac{1}{\tau^2}) \| u^k - p^k \|^2 \\
				\leq & {\rho}^2 \| x^k - y^k \|^2.
			\end{aligned}
		\end{equation}
		From  Theorem \ref{l1}(ii), it follows
		\begin{equation}
			\label{966}
			\| x^k - y^k \| \le\frac{M_1(L+1)}{\varphi_{k-1}+1}+\frac{2L}{\varphi_{k-2} + 1}\|x^0 - x^*\|_{\mathcal{M}},\quad k\ge2.
		\end{equation}
		Putting  \eqref{96} and \eqref{966} together, we have
		\begin{equation*}\label{92}
			\|R(y^k)\| \leq \frac{ M_1 \rho (L+1)}{\varphi_{k-1}+1} + \frac{2\rho L}{\varphi_{k-2} + 1} \|x^0 - x^*\|_{\mathcal{M}}, \quad k \ge 2.
		\end{equation*}
		
		Now, we estimate the complexity result about the primal-dual gap $ G(p^k, q^k)$. Firstly, by the arithmetic-mean inequality, we get
		\begin{equation}\label{4.3}
			\begin{aligned}
				\left\| \frac{u^{k} - p^{k}}{\tau} + \mathbf{K}^*(q^k - v^k) \right\|
				&\leq \sqrt{ \frac{2}{\tau^2} \left\| u^{k} - p^{k} \right\|^2 + 2 \left\| \mathbf{K} \right\|^2 \left\| q^k - v^k \right\|^2 }
				&\leq \rho \left\| x^{k} - y^{k} \right\|,
			\end{aligned}
		\end{equation}
		and
		\begin{equation}\label{4.4}  
			\begin{aligned}
				\left\| \frac{v^{k} - q^{k}}{\sigma} +  \mathbf{K}(p^k -  u^{k}) \right\|
				&\leq\sqrt{2\frac{1}{\sigma^2}\|v^{k} - q^{k}\|^2+2\|\mathbf{K}\|^2\|p^k-u^k\|^2}
				&\leq \rho \|x^{k} - y^{k}\|.
			\end{aligned}
		\end{equation}
		From Lemma \ref{lem41}, \eqref{4.3}, \eqref{4.4} and Cauchy--Schwarz inequality, we have
		\begin{equation*}\label{4.2}
			\begin{aligned}
				G(p^k, q^k)
				&\leq \left\langle p^k - u^*, \frac{u^{k} - p^{k}}{\tau} + \mathbf{K}^*(q^k - v^k) \right\rangle
				+ \left\langle q^k - v^*, \frac{v^{k} - q^{k}}{\sigma} + \mathbf{K}(p^k - u^{k}) \right\rangle \\
				&\leq \left\| p^k - u^* \right\| \left\| \frac{u^{k} - p^{k}}{\tau} + \mathbf{K}^*(q^k - v^k) \right\|
				+ \left\| q^k - v^* \right\| \left\| \frac{v^{k} - q^{k}}{\sigma} + \mathbf{K}(p^k - u^{k}) \right\| \\
				&\leq  M_2\rho \left\| x^k - y^k \right\|,
			\end{aligned}
		\end{equation*}
		which  with  (\ref{966}), yields
		\[
		G(p^k, q^k)\leq\frac{M_1M_2\rho(L+1) }{\varphi_{k-1}+1}+\frac{2M_2 \rho L}{\varphi_{k-2} + 1}\|x^0 - x^*\|_{\mathcal{M}},
		\]
		The proof is completed.
	\end{proof}

\begin{remark}
Recall that Douglas--Rachford (DR) algorithm can be  regarded as a special case of CP algorithm (see, e.g., \cite{CP}).
Consequently,  an accelerated version of DR algorithm can be derived from Algorithm \ref{Alg81}. Correspondingly,  the   convergence rate results established for the accelerated CP algorithm  are readily applicable to its DR counterpart.
\end{remark}
	
	\section{Numerical experiments}\label{sect5}
	To show the performance  of the proposed aCP algorithm and Restarted aCP algorithm, we conduct  numerical tests on two benchmark problems: the minimax matrix game\cite{vonNeumann1928} and the LASSO problem\cite{lasso}. Comparisons are made against the CP algorithm, HCP algorithm and  Restarted HCP algorithm in \cite{Zhang}  under identical experimental setups. All experiments were implemented in Python 3.8 running on a 64-bit Windows PC with an Intel(R)Core(TM)i5-8265U CPU @ 1.60GHz and 8GB of RAM.
	\subsection*{Problem 5.1 (Minimax matrix game)}The minimax matrix game problem is given by
$$
		\min_{u \in \Delta_{q}} \max_{v \in \Delta_{p}} \langle \mathbf{K}u, v \rangle,
$$
	where \( \mathbf{K} \in \mathbb{R}^{p \times q} \), \( \Delta_{q} = \{ u \in \mathbb{R}^{q} : \sum_{i} u_{i} = 1, \, u \geq 0 \} \) and \( \Delta_{p} = \{ v \in \mathbb{R}^{p} : \sum_{i} v_{i} = 1, \, v \geq 0 \} \) denote the standard unit simplex in \( \mathbb{R}^{q} \) and \( \mathbb{R}^{p} \), respectively.
	Let \(\iota_{\Delta_q}\) and \(\iota_{\Delta_p}\) be the indicator functions of \(\Delta_q\) and \(\Delta_p\) respectively. The problem (\ref{51}) can be rewritten as follows
	\begin{equation}
		\label{51}
	\min_{u \in \mathbb{R}^q} \max_{v \in \mathbb{R}^p} \left\{ \iota_{\Delta_q}(u) + \langle  \mathbf{K}u, v \rangle - \iota_{\Delta_p}(v) \right\}.
	\end{equation}
	Clearly, (\ref{51}) is a special case of the saddle-point problem (\ref{saddle}) with \( f = \iota_{\Delta_{q}} \) and \( g^{*} = \iota_{\Delta_{p}} \).
The matrix \( \mathbf{K} \in \mathbb{R}^{p \times q} \) was generated randomly in the following four different ways with random number generator \( {seed} = 50 \).
\begin{figure}[h]
		\centering
		\begin{minipage}{0.58\textwidth}
			\centering
			\includegraphics[width=\linewidth]{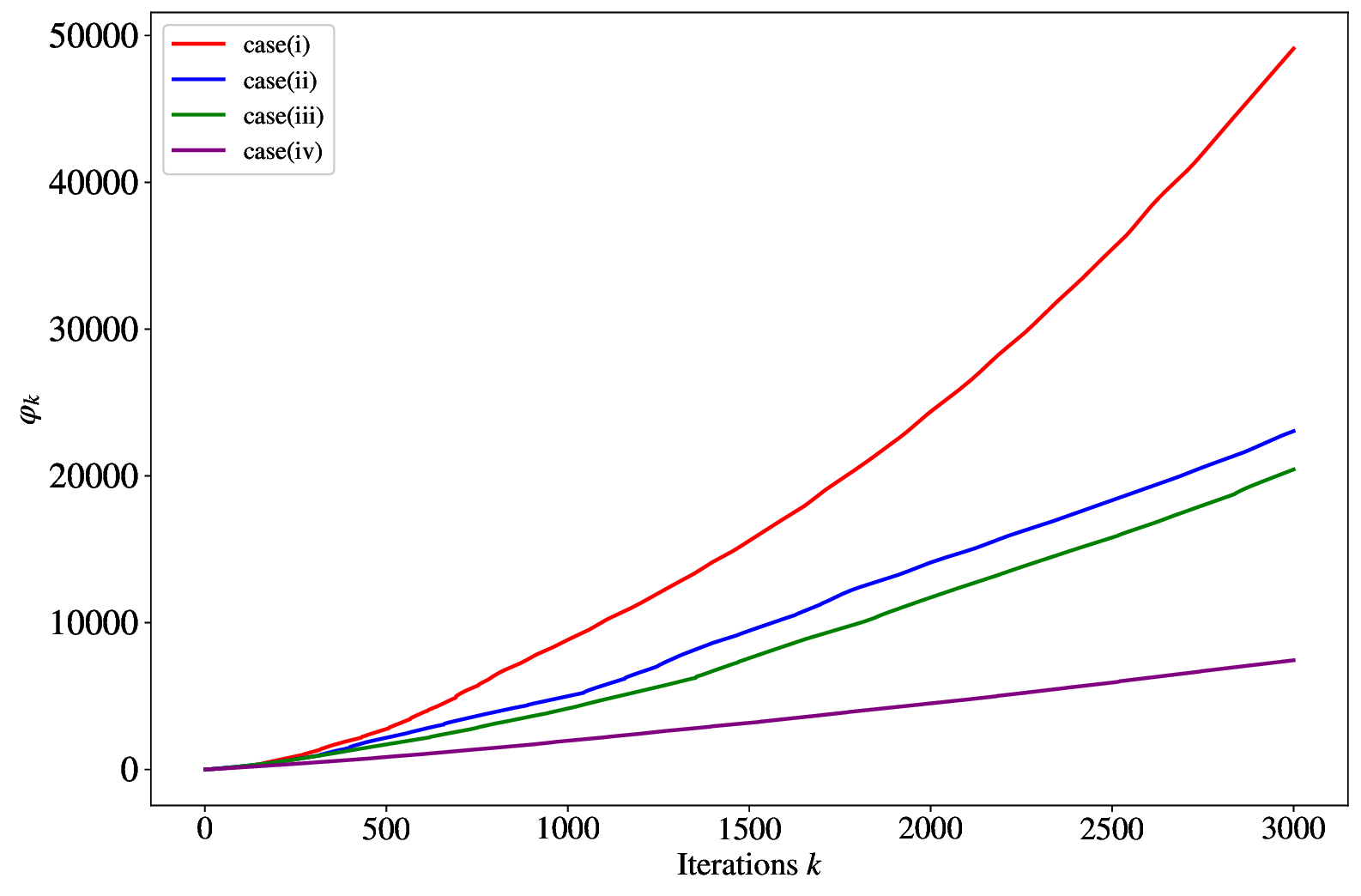}
			\noindent\footnotesize 
			\label{fig:1a}
		\end{minipage}
		\caption{Increase of $\varphi_k$ with $k$ for Problem 5.1}
		\label{fig:1}
	\end{figure}
	\begin{itemize}
		\item[(i)] All entries of \( \mathbf{K} \) were generated independently from the uniform distribution in \([-1, 1]\), and \( (p, q) = (100, 100) \),
		\item[(ii)] All entries of \( \mathbf{K} \) were generated independently from the normal distribution \( \mathcal{N}(0, 1) \) and \( (p, q) = (100, 100) \),
		\item[(iii)] All entries of \( \mathbf{K} \) were generated independently from the normal distribution \( \mathcal{N}(0, 10) \) and \( (p, q) = (500, 500) \),
		\item[(iv)] The matrix \( \mathbf{K} \) is sparse with \( 10\% \) nonzero elements generated independently from the uniform distribution in \([0, 1]\) and \( (p, q) = (1000, 500) \).
	\end{itemize}

	\begin{figure}[t]
		\centering
		\scriptsize  
		
		\begin{subfigure}[b]{0.48\textwidth}
			\centering
			\includegraphics[width=\linewidth]{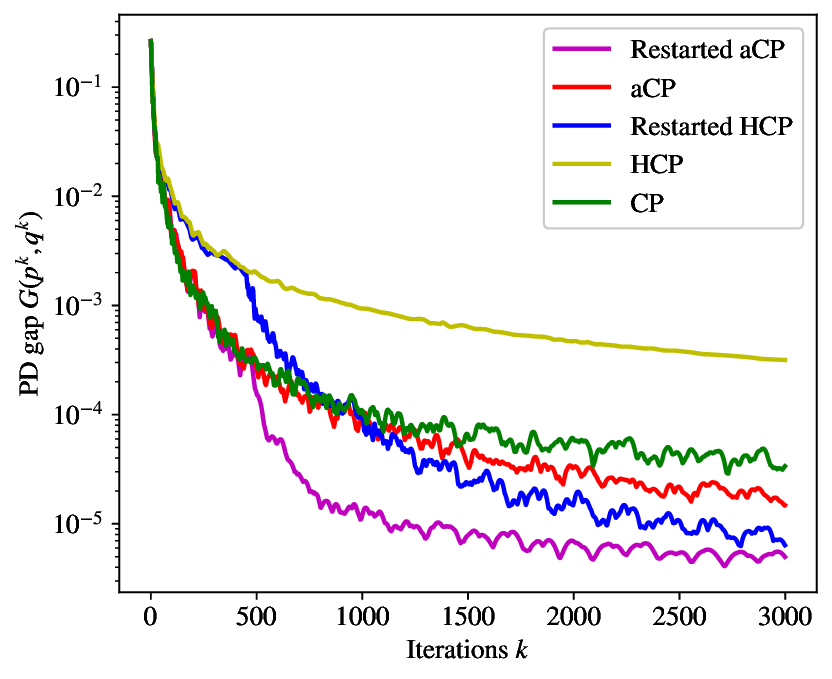}
			\noindent\footnotesize  Test (i) 
			\label{fig:2a}
		\end{subfigure}
		\hfill
		\begin{subfigure}[b]{0.48\textwidth}
			\centering
			\includegraphics[width=\linewidth]{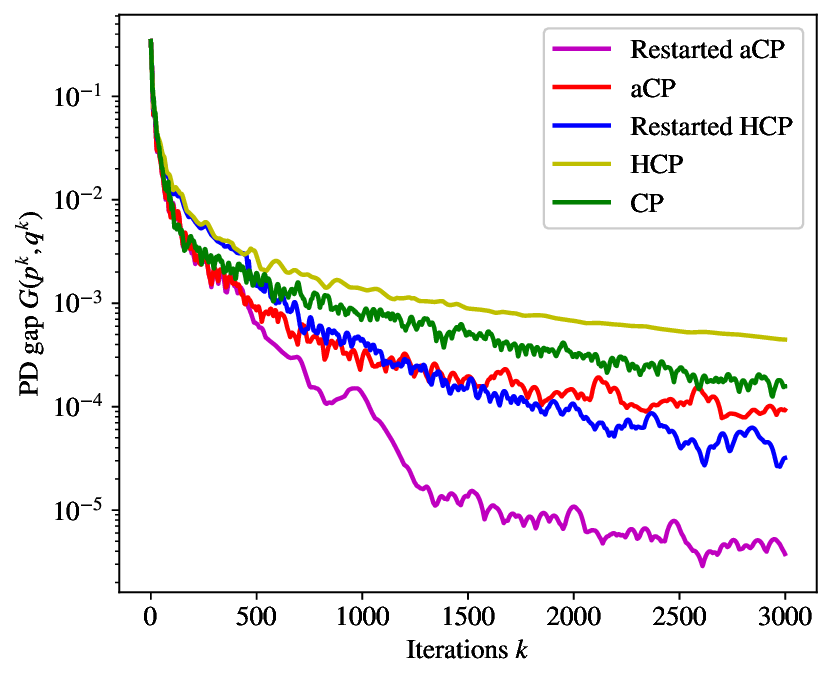}
			\noindent\footnotesize Test (ii) 
			\label{fig:2b}
		\end{subfigure}
		
		\vspace{0mm}  
		\begin{subfigure}[b]{0.48\textwidth}
			\centering
			\includegraphics[width=\linewidth]{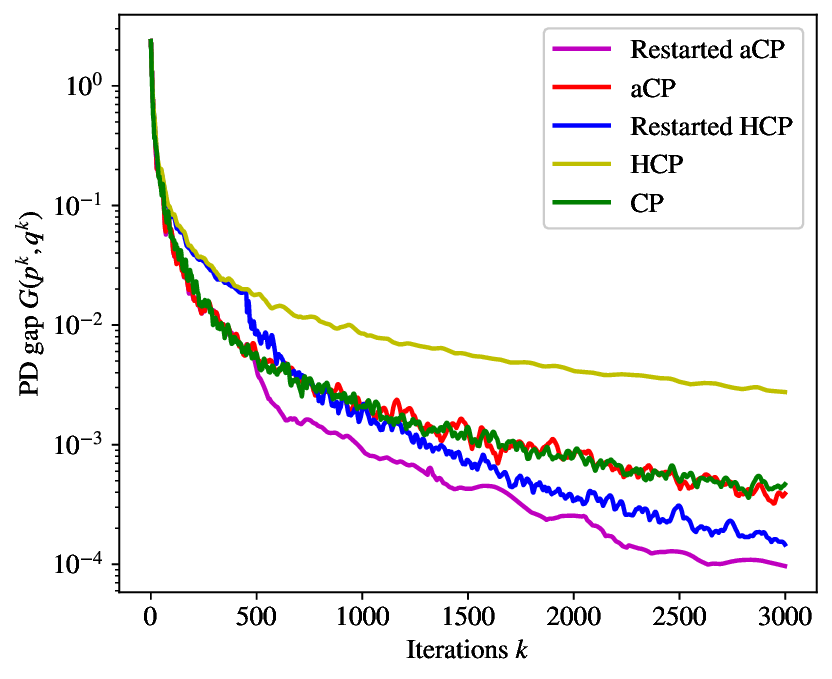}
			\noindent\footnotesize  Test (iii) 
			\label{fig:2c}
		\end{subfigure}
		\hfill
		\begin{subfigure}[b]{0.48\textwidth}
			\centering
			\includegraphics[width=\linewidth]{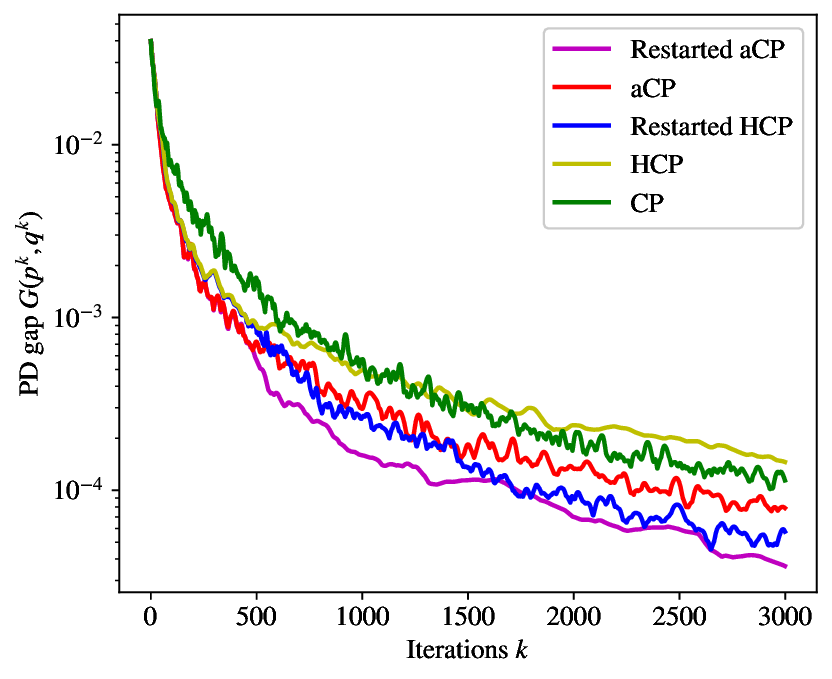}
			\noindent\footnotesize  Test (iv) 
			\label{fig:2d}
		\end{subfigure}
		\vspace{-2mm}
		\caption{Comparison results of all  algorithms on Problem 5.1: PD gap versus iteration numbers}
		\label{fig:2}
	\end{figure}

To evaluate the performance of the different algorithms, we employed the primal-dual gap (PD gap) function  defined in (\ref{pd}), which can be easily computed for a feasible pair  $ (u,v) \in \Delta_{q} \times \Delta_{p}$  by
	\[ G(u, v) := \max_{1 \leq i\leq p} (\mathbf{K}u)_{i} - \min_{1\leq j\leq q} (\mathbf{K}^{*} v)_{j}.
	\]
We also compared the residual   as follows
 $$
 \|R(u,v)\|=\sqrt{\|(u - P_{\Delta_{q}}  (u -  \mathbf{K}^* v)\|^2+\|v - P_{\Delta_{p}} (v +  \mathbf{K} u)\|^2}.
 $$

	Initial points for  all algorithms were set to be \( u_{0} = \frac{1}{q} (1, \ldots, 1)^{\top} \in \mathbb{R}^{q} \) and \( v_{0} = \frac{1}{p} (1, \ldots, 1)^{\top} \in \mathbb{R}^{p} \).  We restarted aCP algorithm and HCP algorithm  every 450
iterations in  the restarted aCP algorithm and  restarted HCP algorithm, respectively. The projection onto the unit simplex was computed by the algorithm from \cite{Duchi2008}.
	We set \( \tau = \sigma = 1 / \| \mathbf{K} \| \) for all algorithms.

\begin{figure}[h]
		\centering
		\begin{subfigure}[b]{0.48\textwidth}
			\centering
			\includegraphics[width=\linewidth]{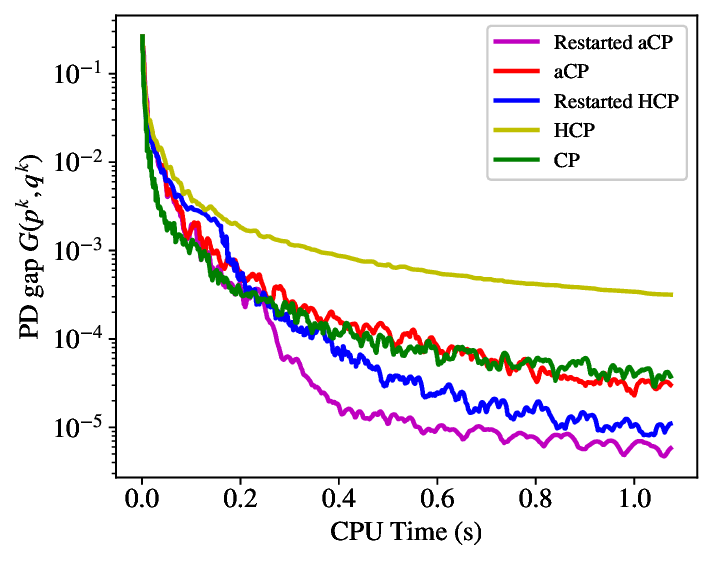}
			\noindent\footnotesize Test (i) 
			\label{fig:3a}
		\end{subfigure}
		\hfill
		\begin{subfigure}[b]{0.48\textwidth}
			\centering
			\includegraphics[width=\linewidth]{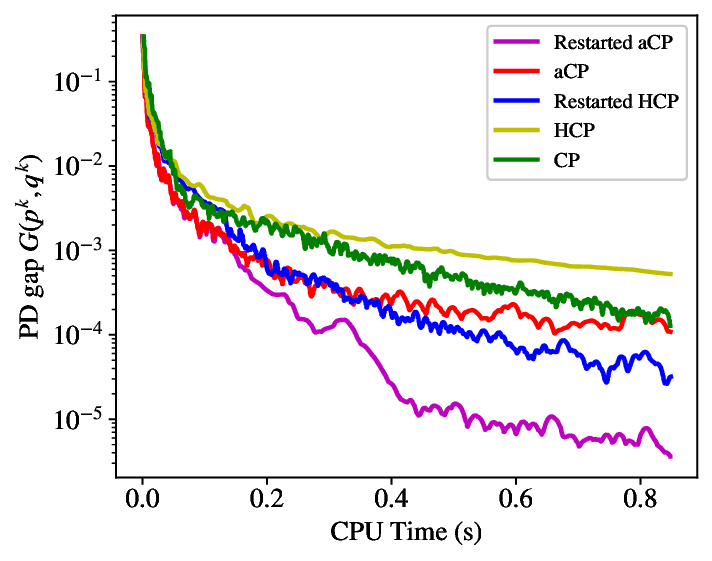}
			\noindent\footnotesize  Test (ii) 
			\label{fig:3b}
		\end{subfigure}
		\vspace{0mm}  
		\begin{subfigure}[b]{0.48\textwidth}
			\centering
			\includegraphics[width=\linewidth]{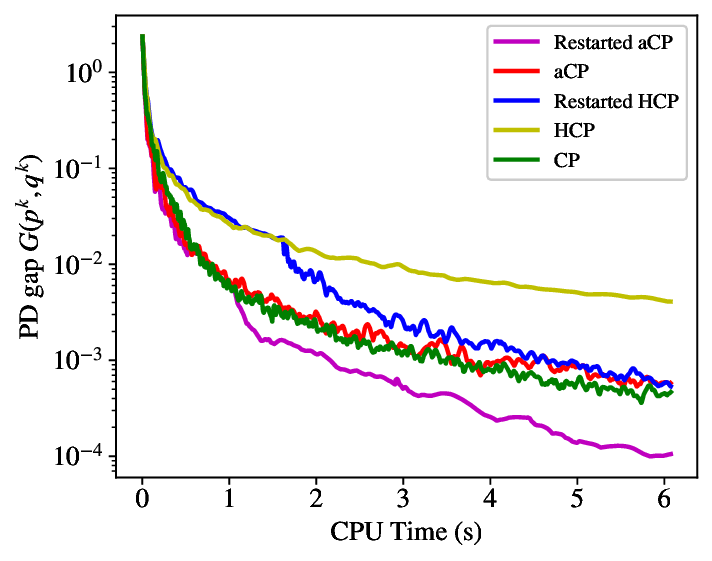}
			\noindent\footnotesize  Test (iii) 
			\label{fig:3c}
		\end{subfigure}
		\hfill
		\begin{subfigure}[b]{0.48\textwidth}
			\centering
			\includegraphics[width=\linewidth]{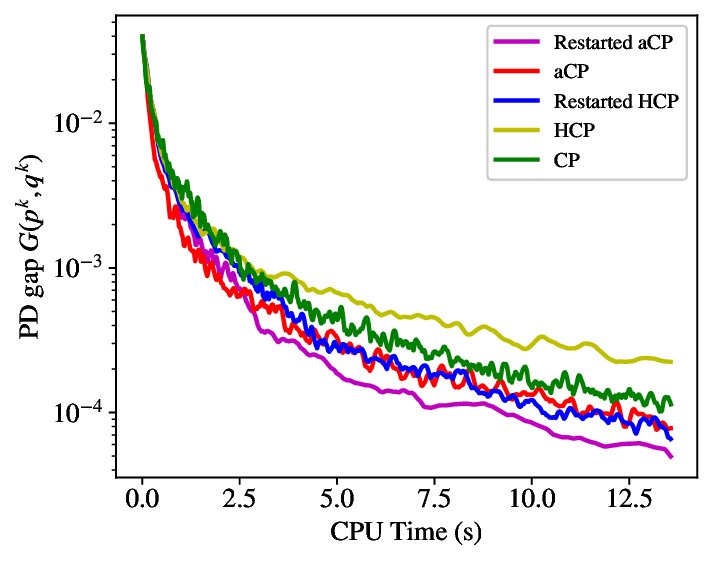}
			\noindent\footnotesize  Test (iv) 
			\label{fig:3d}
		\end{subfigure}
		\vspace{-2mm}
		\caption{Comparison results of all algorithms  on Problem 5.1: PD gap versus CPU times}
		\label{fig:3}
	\end{figure}

\begin{figure}[h]
		\centering
		\scriptsize  
		
		\begin{subfigure}[b]{0.48\textwidth}
			\centering
			\includegraphics[width=\linewidth]{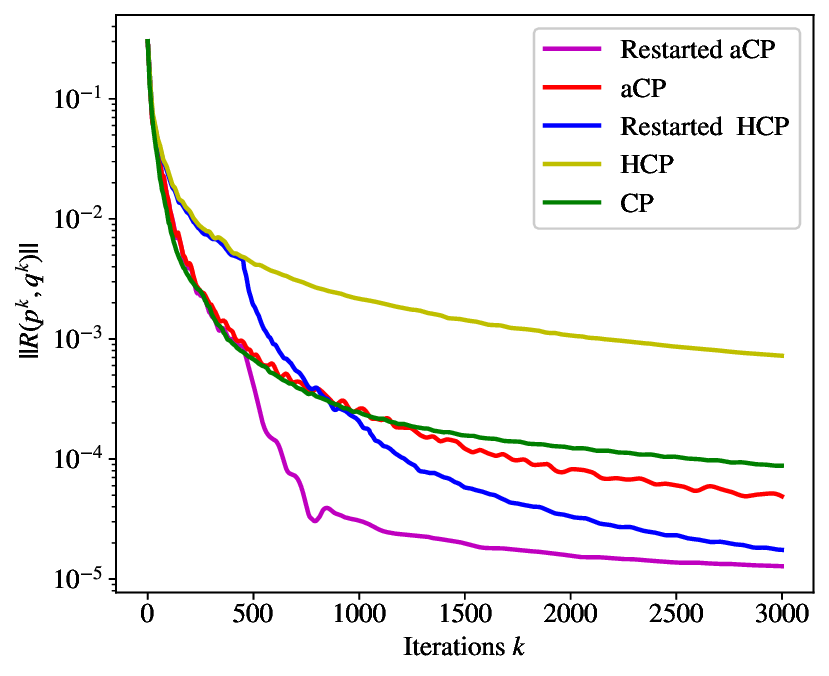}
			\noindent\footnotesize  Test (i) 
			\label{fig:2a}
		\end{subfigure}
		\hfill
		\begin{subfigure}[b]{0.48\textwidth}
			\centering
			\includegraphics[width=\linewidth]{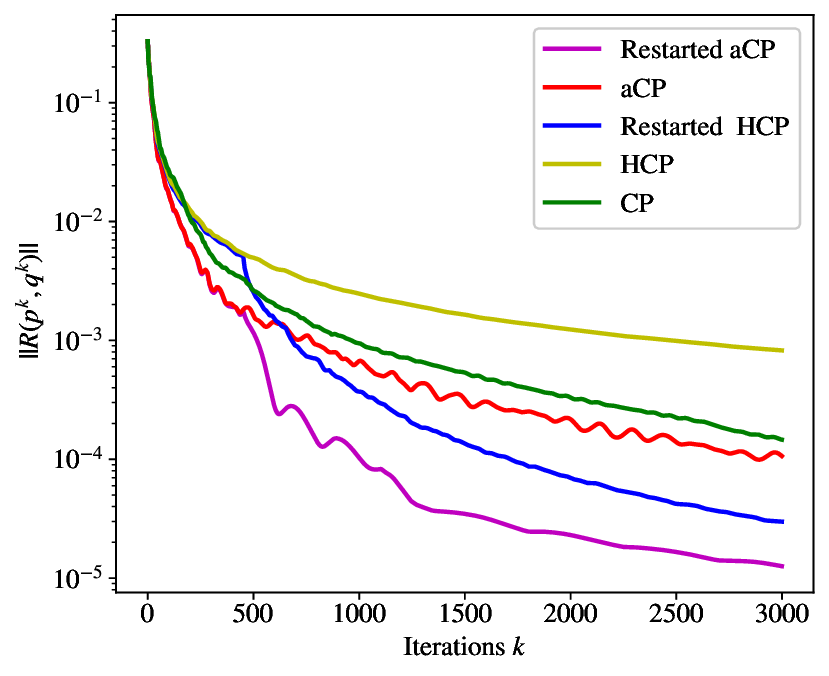}
			\noindent\footnotesize Test (ii) 
			\label{fig:2b}
		\end{subfigure}
		
		\vspace{0mm}  
		\begin{subfigure}[b]{0.48\textwidth}
			\centering
			\includegraphics[width=\linewidth]{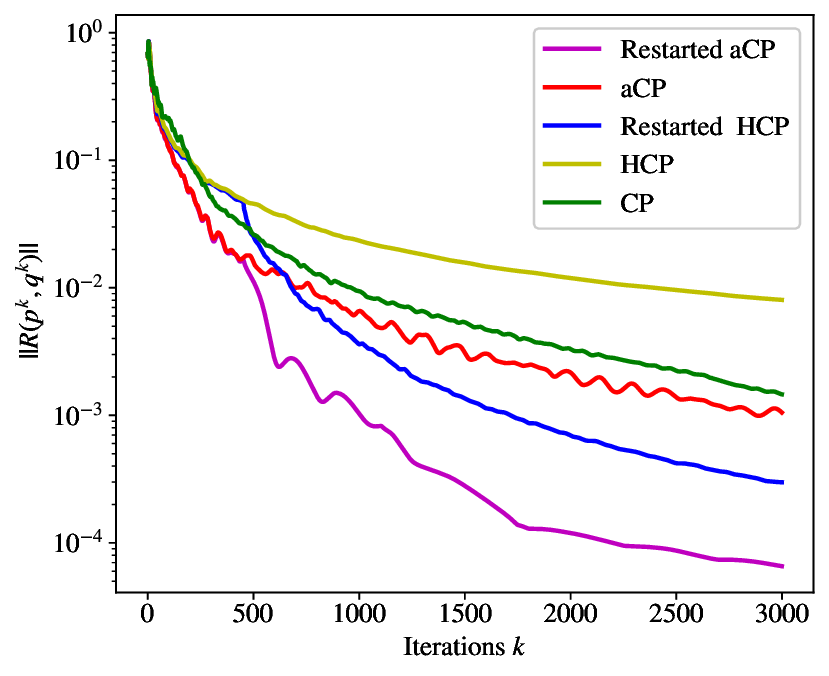}
			\noindent\footnotesize  Test (iii) 
			\label{fig:2c}
		\end{subfigure}
		\hfill
		\begin{subfigure}[b]{0.48\textwidth}
			\centering
			\includegraphics[width=\linewidth]{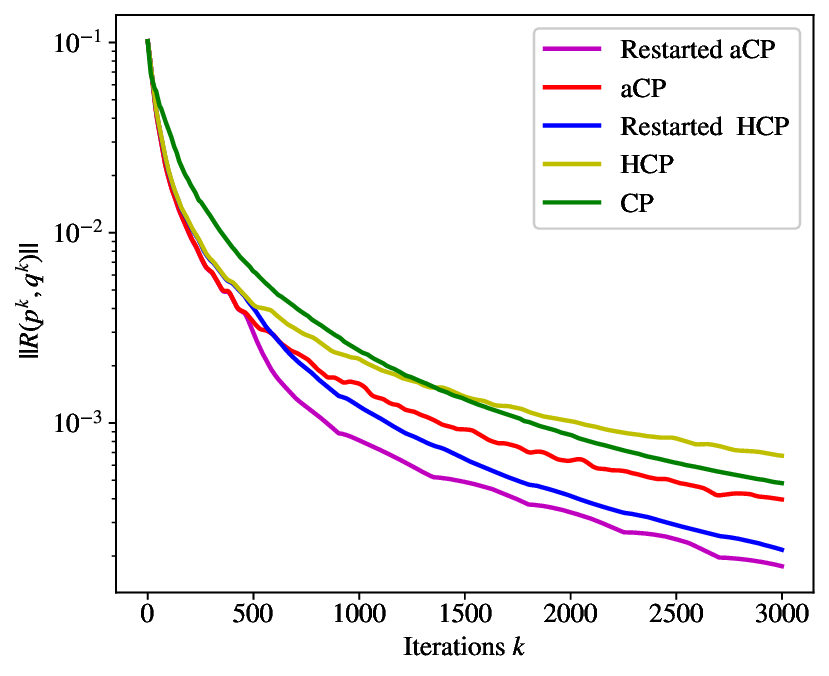}
			\noindent\footnotesize  Test (iv) 
			\label{fig:2d}
		\end{subfigure}
		\vspace{-2mm}
		\caption{Comparison results of all  algorithms on Problem 5.1: Residual   versus iteration numbers}
		\label{fig:4}
	\end{figure}

\begin{figure}[h]
		\centering
		\begin{subfigure}[b]{0.48\textwidth}
			\centering
			\includegraphics[width=\linewidth]{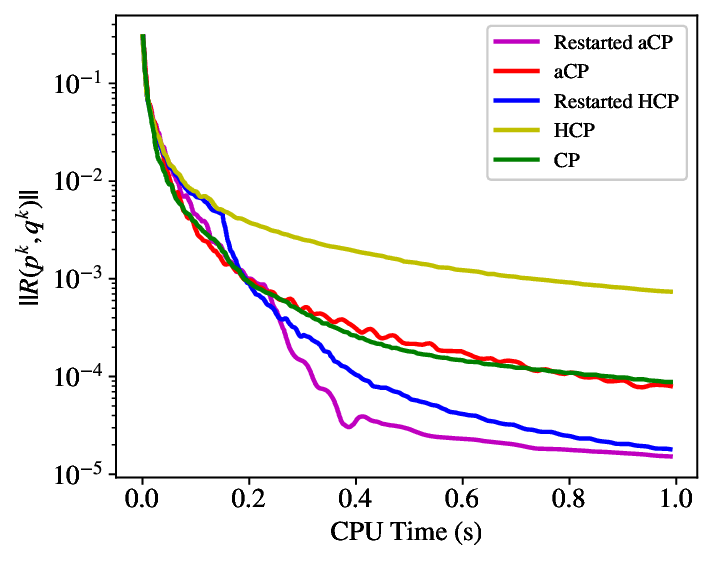}
			\noindent\footnotesize Test (i) 
			\label{fig:3a}
		\end{subfigure}
		\hfill
		\begin{subfigure}[b]{0.48\textwidth}
			\centering
			\includegraphics[width=\linewidth]{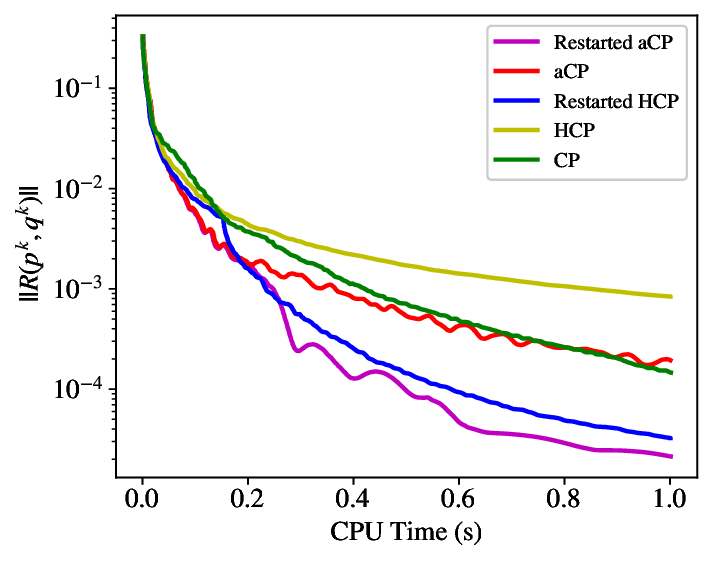}
			\noindent\footnotesize  Test (ii) 
			\label{fig:3b}
		\end{subfigure}
		\vspace{0mm}  
		\begin{subfigure}[b]{0.48\textwidth}
			\centering
			\includegraphics[width=\linewidth]{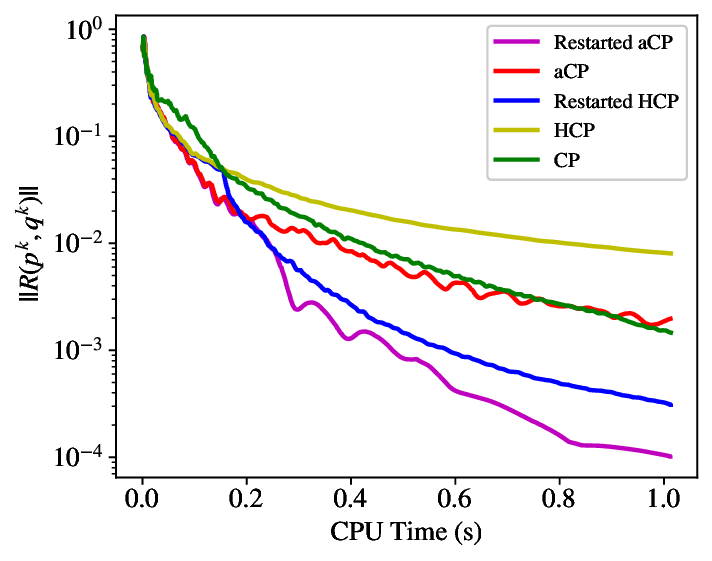}
			\noindent\footnotesize  Test (iii) 
			\label{fig:3c}
		\end{subfigure}
		\hfill
		\begin{subfigure}[b]{0.48\textwidth}
			\centering
			\includegraphics[width=\linewidth]{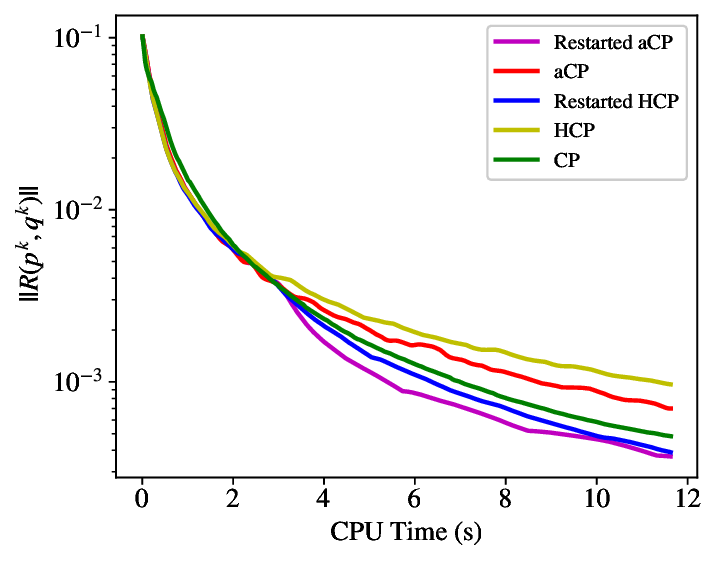}
			\noindent\footnotesize  Test (iv) 
			\label{fig:3d}
		\end{subfigure}
		\vspace{-2mm}
		\caption{Comparison results of all algorithms  on Problem 5.1: Residual   versus CPU times}
		\label{fig:42}
	\end{figure}

	Figure \ref{fig:1} exhibits increase of the anchoring parameter $\varphi_k$  with $k$. The decreasing behavior of the PD gap function as the algorithms progress is presented in both Figures \ref{fig:2} and \ref{fig:3}.  It can be observed from these two figures that the Restarted aCP algorithm behaves better than other algorithms  both in terms of iterations and CPU time. Figures \ref{fig:4} and \ref{fig:42} show that the residual  declines  as each algorithm proceeds, and the restarted aCP algorithm clearly  outperforms  the other algorithms.

	\subsection*{Problem 5.2 (LASSO)}
	The LASSO  problem  has the form
	\begin{equation*}
		\label{422}
		\min_u F(u) := \frac{1}{2} \| \mathbf{K}u - b \|^2 + \mu \| u \|_1,
	\end{equation*}
	where  \( \mathbf{K} \in \mathbb{R}^{p \times q} \) and  \( b \in \mathbb{R}^p \) are known,   \( u \in \mathbb{R}^q \) is an unknown signal, and $p$ is much smaller than $q$ in compressive sensing applications.
\begin{figure}[h]
		\centering
		\begin{minipage}{0.58\textwidth}
			\centering
			\includegraphics[width=\linewidth]{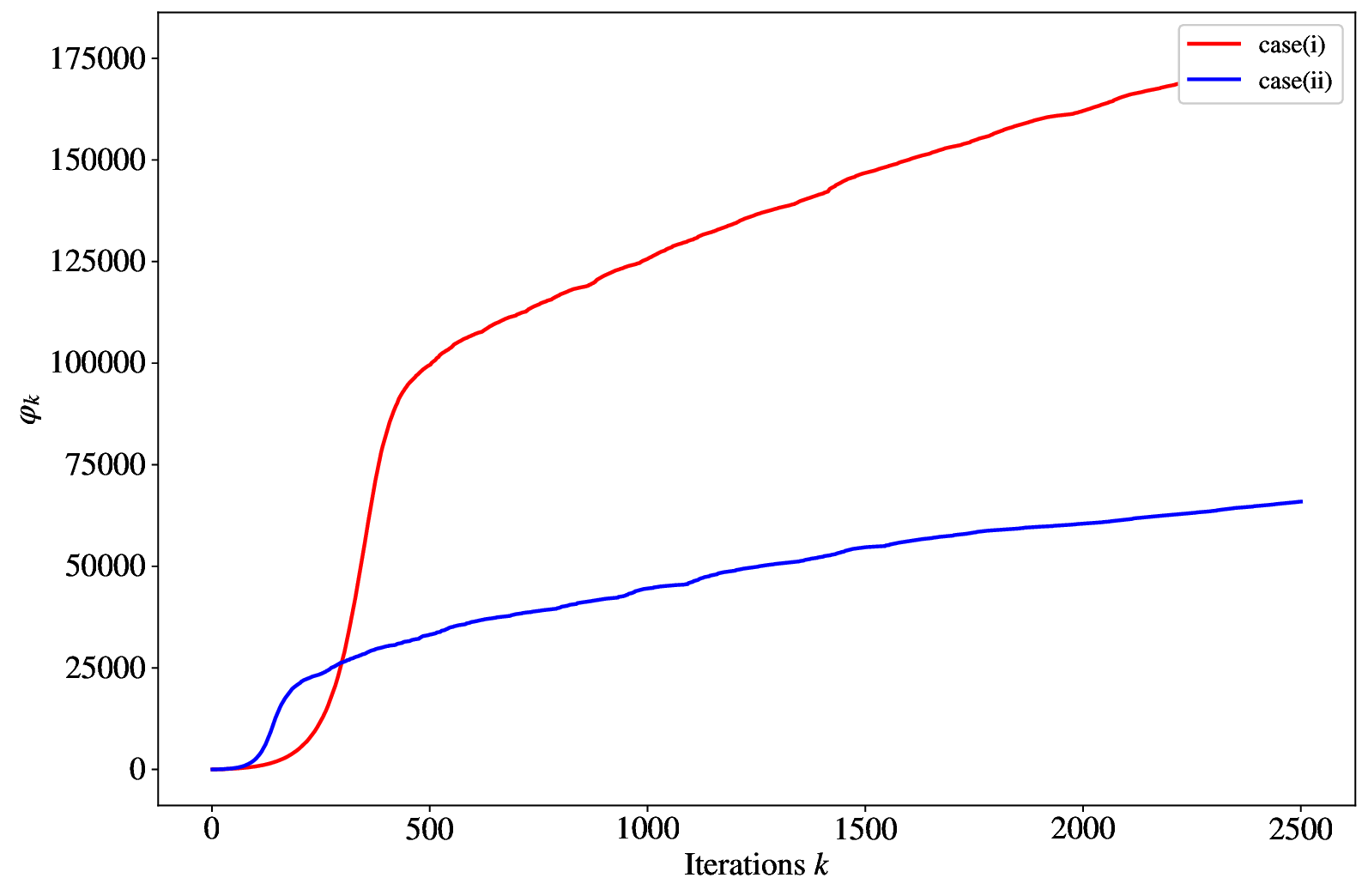}
			\noindent\footnotesize 
			\label{fig:1a}
		\end{minipage}
		\caption{Increase of $\varphi_k$ with $k$ for Problem 5.2}
		\label{fig:5}
\end{figure}
	The above LASSO problem  corresponds to \eqref{tv1}
	with $ f(\cdot) = \mu \| \cdot \|_1$ and  $g(\cdot) = \frac{1}{2} \| \cdot - b \|_2^2 $.
	
	\begin{figure}[h]
		\centering
		\begin{subfigure}[b]{0.48\textwidth}
			\centering
			\includegraphics[width=\textwidth]{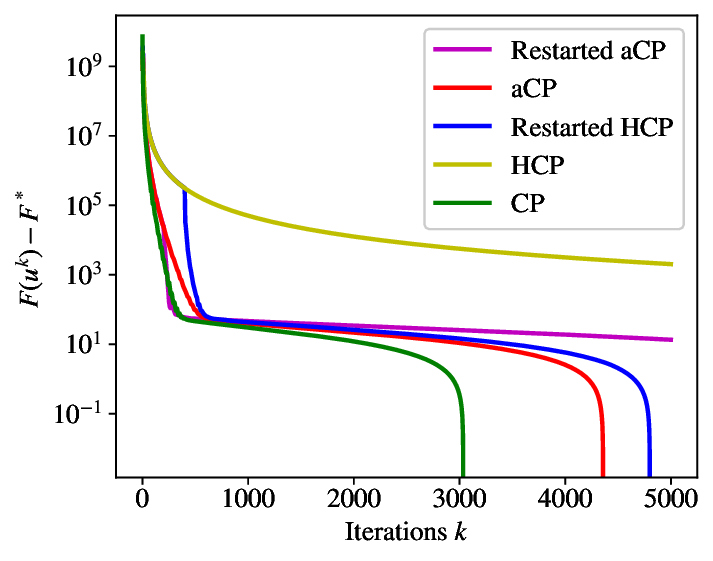}
			\noindent\footnotesize Test (i) 
			\label{fig:5a}
		\end{subfigure}
		\hfill
		\begin{subfigure}[b]{0.48\textwidth}
			\centering
			\includegraphics[width=\textwidth]{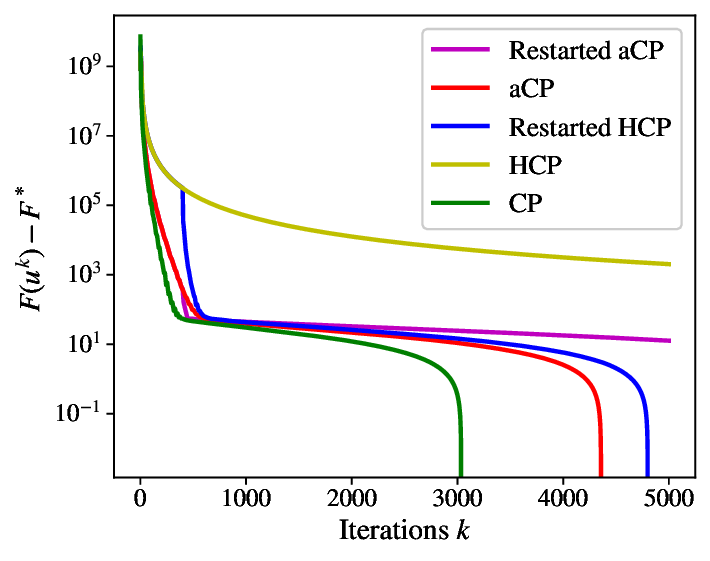}
			\noindent\footnotesize Test (ii) 
			\label{fig:5b}
		\end{subfigure}
		\vspace{-2mm}
		\caption{Comparison of the decay of function errors with iteration numbers among all  algorithms  for LASSO}
		\label{fig:22}
	\end{figure}
	We set $ seed $ = 100 and generated an exact solution \( u^* \in \mathbb{R}^q \) randomly. Specifically, $ s $ nonzero components of \( u^* \) were determined uniformly at random, and their values were drawn from the uniform distribution in \([-10, 10]\).  We generated additive noise \( \nu \in \mathbb{R}^p \)  from \( \mathcal{N}(0, 0.1) \) and set \( b = \mathbf{K}u^* + \nu \). The matrix \( \mathbf{K} \in \mathbb{R}^{p \times q} \) was generated as in \cite{Malitsky}.
	\begin{enumerate}[label=(\roman*)]
		\item All entries of \( \mathbf{K} \) were generated independently from \( \mathcal{N}(0, 1) \), the normal distribution with mean 0 and standard deviation 1.
		\item First, we generated a matrix \( A \in \mathbb{R}^{p \times q} \), whose entries are independently drawn from \( \mathcal{N}(0, 1) \). Then, for a scalar \( v \in (0, 1) \), we constructed the matrix \( \mathbf{K} \) column by column as follows: \( \mathbf{K}_1 = \frac{A_1}{\sqrt{1 - v^2}} \) and \( \mathbf{K}_j = v \mathbf{K}_{j-1} + A_j \), \( j = 2, \ldots, q \). Here \( \mathbf{K}_j \) and \( A_j \) represent the \( j \)-th column of \( \mathbf{K} \) and \( A \), respectively. As \( v \) increases, \( \mathbf{K} \) becomes more ill-conditioned. In this experiment, we tested \( v = 0.5 \).
	\end{enumerate}
In both cases, the regularization parameter $\mu$ was set to be 0.1 and ($q$, $p$, $s$)=(2000, 1000, 100).

All algorithms adopt the initial points  $u^0=0$ and $v^0=-b$ which are a common choice since they satisfy the constraint relation of the corresponding primal-dual problem and  are consistent with the strong sparsity of the problem’s exact solution.  We restarted aCP algorithm and HCP algorithm  every 400
iterations in  the restarted aCP algorithm and  restarted HCP algorithm, respectively.

	\begin{figure}[h]
		\begin{subfigure}[b]{0.48\textwidth}
			\centering
			\includegraphics[width=\textwidth]{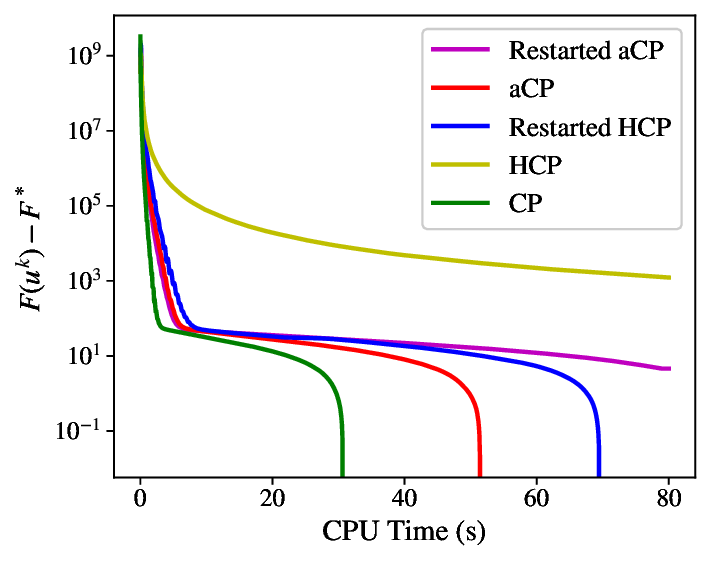}
			\noindent\footnotesize Test (i) 
			\label{fig:6a}
		\end{subfigure}
		\hfill
		\begin{subfigure}[b]{0.48\textwidth}
			\centering
			\includegraphics[width=\textwidth]{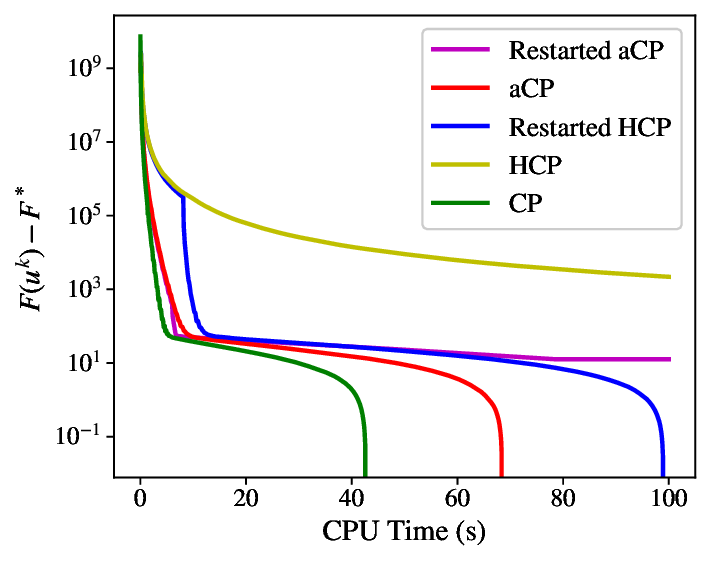}
			\noindent\footnotesize Test (ii) 
			\label{fig:6b}
		\end{subfigure}
		\vspace{-1mm}
		\caption{Comparison of the decay of function errors with the CPU time results of all  algorithms for LASSO}
		\label{fig:23}
	\end{figure}
	\begin{figure}[h]
		\centering
		\begin{subfigure}[b]{0.48\textwidth}
			\centering
			\includegraphics[width=\textwidth]{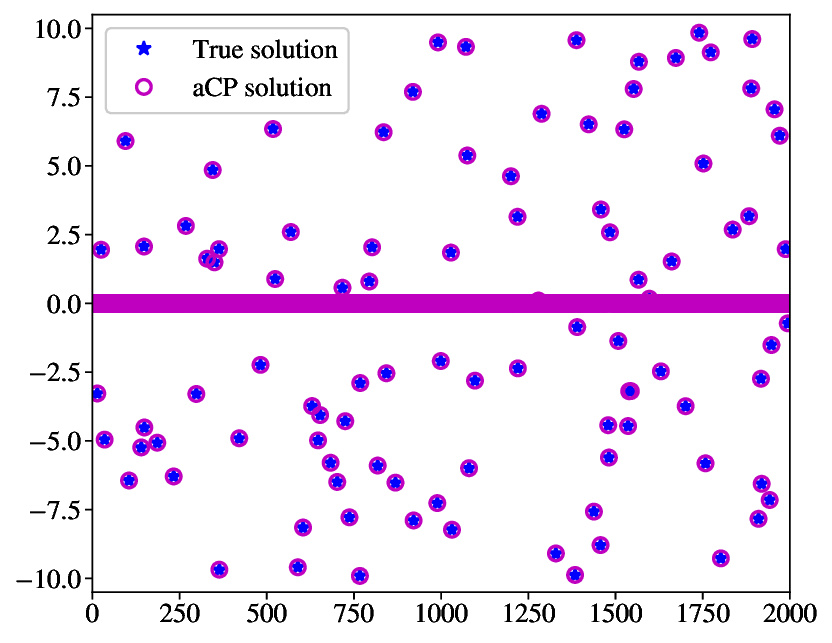}
			\noindent\footnotesize Test (i) 
			\label{fig:7a}
		\end{subfigure}
		\hfill
		\begin{subfigure}[b]{0.48\textwidth}
			\centering
			\includegraphics[width=\textwidth]{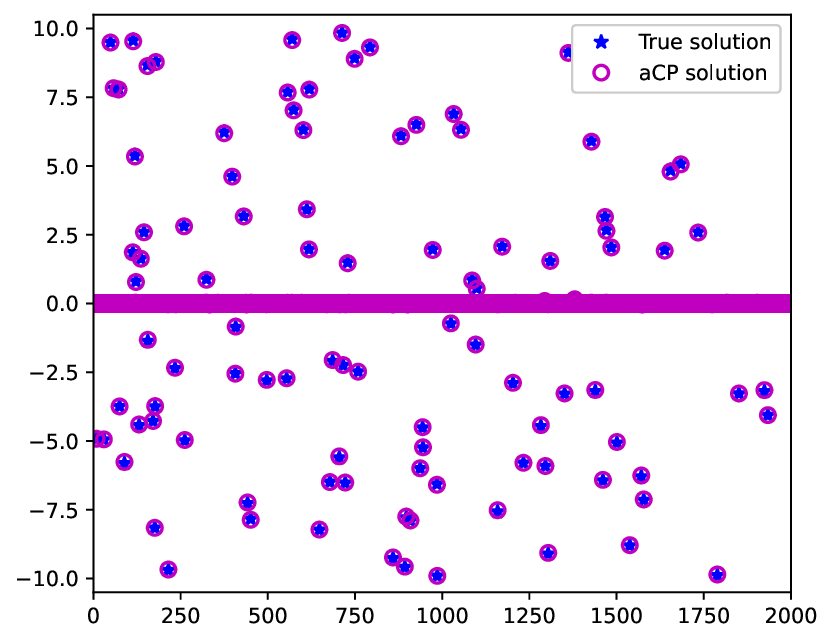}
			\noindent\footnotesize Test (ii) 
			\label{fig:7b}
		\end{subfigure}
		\vspace{-1mm}
		\caption{The true solutions and the solutions by the aCP for LASSO}
		\label{fig:33}
	\end{figure}Figure \ref{fig:5} shows the increase of the anchoring parameter $\varphi_k$  with $k$.  The decay of the errors $F(u^k) - F^*$ of the function values  as the algorithms proceed is presented in Figures \ref{fig:22} and  \ref{fig:23}, where  $F^* = F(u^*)$ is the value of objective function $F$ at the exact solution $ u^* $.  It can be observed that the CP algorithm performs best, with the aCP algorithm ranking second. Moreover, the CP algorithm and the aCP algorithm outperform the HCP algorithm in the first few iterations.  The excellent performance of the CP algorithm and the restart-free property of the aCP algorithm can likely be attributed to the appropriate choice of initial points. Notably, in this experiment, the restarted HCP algorithm always  improves the original HCP algorithm,  which is consistent with the numerical results in \cite{Zhang}.

	To illustrate the ability of recovering the original sparse solution by the aCP algorithm, we plot in Figure~\ref{fig:33} the true solutions and the approximate solutions on a random instance $(q,p,s)=(2000,1000,100)$. The true solution is represented by asterisks, while circles are the estimates obtained by the aCP algorithm. We see that the estimates obtained by the aCP algorithm  are quite close to the true values.

\vspace{5cm}
	\section{Conclusion}
	
	In this article, we study a  preconditioned Halpern iteration with adaptive anchoring parameters and its acceleration to the degenerate PPP method. We establish the strong convergence of the proposed methods and present  convergence rate of at least $\mathcal{O}(1/k)$. Moreover, we focus on an accelerated CP algorithm achieving at least $\mathcal{O}(1/k)$ convergence rate
concerning the residual mapping and the primal-dual gap.
We  consider the acceleration of the degenerate PPP method as it encompasses many splitting methods and their variants besides the Chambolle--Pock algorithm.
This article focuses on investigating the convergence rates for the accelerated CP algorithm, with the exploration of convergence rates for other degenerate PPP  methods deferred to future work.
\\

\section*{Acknowledgements}

We are deeply grateful to the anonymous referee for the constructive suggestions and critical comments on the manuscript, especially for supplementing the argument for the existence of \( P_{\operatorname{Fix}(T)}^{\mathcal{M}}(u) \). These  suggestions and comments have significantly helped us improve the presentation of the paper. We are grateful to Professor Songnian He for valuable discussions.\\

\noindent$\textbf{Funding}$ The second author was supported by National Natural Science Foundation of China (No.12271273).\\

\noindent{$\textbf{Data Availability}$ Not applicable.}\\

\subsection*{Declarations}

\noindent{$\textbf{Ethics approval and consent to participate}$ Not applicable.}\\	

\noindent{$\textbf{Consent for publication}$ Not applicable.}\\

\noindent{$\textbf{Conflict~of~Interest}$ The authors declare that they have no conflict of interest.}


\begin{thebibliography}{99}
		
		\bibitem{Bertsekas1982}
		D. P. Bertsekas, E. M. Gafni.
		\textit{Projection methods for variational inequalities with application to the traffic assignment problem}.
		\emph{Mathematical Programming Studies}, \textbf{17}, 1982, pp. 139--159.
		
		
		
		\bibitem{Bouwmans2016}
		T. Bouwmans, N. S. Aybat, E. H. Zahzah.
		\textit{Handbook of robust low-rank and sparse matrix decomposition: Applications in image and video processing}.
		\emph{ CRC Press, Taylor and Francis Group}, \textbf{05 }, 2016.
		45, 2016.
		
		\bibitem{Hayden2013}
		S. Hayden, O. Stanley.
		\textit{A low patch-rank interpretation of texture}.
		\emph{SIAM Journal on Imaging Sciences}, \textbf{6}(1), 2013, pp. 226--262.

\bibitem{CP}
		A. Chambolle, T. Pock.
		\textit{A first-order primal-dual algorithm for convex problems with applications to imaging}.
		\emph{Journal of Mathematical Imaging and Vision}, \textbf{40}(1), 2011, pp. 120--145.



		\bibitem{Gabay1976}
		D. Gabay, B. Mercier.
		\textit{A dual algorithm for the solution of nonlinear variational problems via finite-element approximations}.
		\emph{Computers \& Mathematics with Applications}, \textbf{2}, 1976, pp. 17--40.
		
		
		\bibitem{Jonathan1992}
		J. Eckstein, D. P. Bertsekas.
		\textit{On the Douglas--Rachford splitting method and the proximal point algorithm for maximal monotone operators}.
		\emph{Mathematical Programming}, \textbf{55}(1-3), 1992, pp. 293--318.
		
		\bibitem{Liu2021}
		Y. Liu, Y. Xu, W. Yin.
		\textit{Acceleration of primal-dual methods by preconditioning and simple subproblem procedures}.
		\emph{Journal of Scientific Computing}, \textbf{86}(2), 2021.
		

		\bibitem{Malitsky}
		Y. Malitsky, T. Pock.
		\textit{A first-order primal-dual algorithm with linesearch}.
		\emph{SIAM Journal on Optimization}, \textbf{28}(1), 2018, pp. 411--432.
		
\bibitem{Bauschke2017}
		H. H. Bauschke, P. L. Combettes.
		\textit{Convex analysis and monotone operator theory in Hilbert spaces}.
		\emph{Springer}, 2017.
		
		\bibitem{Bredies}
		K. Bredies, E. Chenchene, D. A. Lorenz, E. Naldi.
		\textit{Degenerate preconditioned proximal point algorithms}.
		\emph{SIAM Journal on Optimization}, \textbf{32}(4), 2022, pp. 2376--2401.


		
		\bibitem{Douglas1956}
		J. Douglas, H. H. Rachford.
		\textit{On the numerical solution of heat conduction problems in two and three space variables}.
		\emph{Transactions of the American Mathematical Society}, \textbf{82}(2), 1956, pp. 421--439.
		
		\bibitem{Peaceman1955}
		D. W. Peaceman, H. H. Rachford.
		\textit{The numerical solution of parabolic and elliptic differential equations}.
		\emph{Journal of the Society for Industrial and Applied Mathematics}, \textbf{3}(1), 1955, pp. 28--41.
		
		\bibitem{Davis2017}
		D. Davis, W. Yin.
		\textit{A three-operator splitting scheme and its optimization applications}.
		\emph{Set-Valued and Variational Analysis}, \textbf{25}(4), 2017, pp. 829--858.
		
		
		\bibitem{Yang2025}
		B. Yang, X. Zhao, X. Li, D.  Sun.
		\textit{An accelerated proximal alternating direction method of multipliers for optimal decentralized control of uncertain systems}.
		\emph{Journal of Optimization Theory and Applications}, \textbf{204}(9), 2025

\bibitem{Aragon2023}
F. J. Aragón-Artacho, Y. Malitsky, M. K. Tam, D. Torregrosa-Belén.
\textit{Distributed forward-backward methods for ring networks}.
\emph{Computational Optimization and Applications}, \textbf{86}(3), 2023, pp. 845--870.


       \bibitem{Ryu2020}
       E. K. Ryu.
       \textit{Uniqueness of DRS as the 2-operator resolvent-splitting and impossibility of 3-operator resolvent-splitting}.
       \emph{Mathematical Programming}, \textbf{182}, 2020, pp. 233--273.

    \bibitem{Malitsky2023}
	Y. Malitsky, M. K. Tam.
	\textit{Resolvent splitting for sums of monotone operators with minimal lifting}.
	\emph{Mathematical Programming}, \textbf{201}, 2023, pp. 231--262.


       \bibitem{Campoy2022}
	R. Campoy.
	\textit{A product space reformulation with reduced dimension for splitting algorithms}.
	\emph{Computational Optimization and Applications}, \textbf{83}, 2022, pp. 319--348.




\bibitem{Bredies2024}
	K. Bredies, E. Chenchene, E. Naldi.
	\textit{Graph and distributed extensions of the Douglas--Rachford method}.
	\emph{SIAM Journal on Optimization}, \textbf{34}(2), 2024, pp. 1569--1594.

\bibitem{Akerman2025}
      A. {\AA}kerman, E. Chenchene, P. Giselsson, E. Naldi.
     \textit{Splitting the forward-backward algorithm: a full characterization}.
      \emph{arXiv preprint arXiv:2504.10999}, 2025.


\bibitem{Aragón-Artacho2024}
        F. J. Aragón-Artacho, R. Campoy, C. López-Pastor.
       \textit{Forward-backward algorithms devised by graphs}.
      \emph{SIAM Journal on Optimization}, \textbf{35}(4), 2025, pp. 2423--2451.


\bibitem{Lorenz2024}
        D. A. Lorenz, J. Marquardt, E. Naldi.
       \textit{The degenerate variable metric proximal point algorithm and adaptive stepsizes for primal-dual Douglas--Rachford}.
      \emph{Optimization} \textbf{74}(6), 2024, pp. 1355--1381.

     \bibitem{MarquesAlves2024}
        M. M. Alves, D. A. Lorenz, E. Naldi.
       \textit{A general framework for inexact splitting algorithms with relative errors and applications to Chambolle--Pock and Davis--Yin methods}.
      \emph{Computational Optimization and Applications}, \textbf{93}, 2026, pp. 729--763.
	
		
		\bibitem{Sun2025}
		D. Sun, Y. Yuan, G. Zhang, X. Zhao.
		\textit{Accelerating preconditioned ADMM via degenerate proximal point mappings}.
		\emph{SIAM Journal on Optimization}, \textbf{35}, 2025, pp. 1165--1193.
		
		\bibitem{Halpern1967}
		B. Halpern.
		\textit{Fixed points of nonexpanding maps}.
		\emph{Bulletin of the American Mathematical Society}, \textbf{73}(6), 1967, pp. 957--961.

		\bibitem{Bot2023}
		R. I. Bo\c{t}, D. K. Nguyen.
		\textit{Fast Krasnosel'ski\v{\i}--Mann algorithm with a convergence rate of the fixed point iteration of \(o(1/k)\)}.
		\emph{SIAM Journal on Numerical Analysis}, \textbf{61}(6), 2023, pp. 2813--2843.


      \bibitem{Bot2024}
    R. I. Bo\c{t}, E. Chenchene, J. M. Fadili.
    \textit{Generalized fast Krasnosel'ski\v{i}--Mann method with preconditioners}.
    \emph{arXiv preprint arXiv:2411.18574}, 2024.

		
		
		\bibitem{Zhang}
		S. Zhang, H. Zhang, H. Wang.
		\textit{HPPP: Halpern-type preconditioned proximal point algorithms and applications to image restoration}.
		\emph{SIAM Journal on Imaging Sciences}, \textbf{18}(2), 2025, pp. 1493--1521.
		
		
		\bibitem{Rudin1992}
		L. I. Rudin, S. Osher, E. Fatemi.
		\textit{Nonlinear total variation based noise removal algorithms}.
		\emph{Physica D: Nonlinear Phenomena}, \textbf{60}, 1992, pp. 259--268.
		
		
		\bibitem{He}
		S. He, H. K. Xu, Q. L. Dong, N. Mei.
		\textit{Convergence analysis of the Halpern iteration with adaptive anchoring parameters}.
		\emph{Mathematics of Computation}, \textbf{93}(345), 2024, pp. 327--345.
		
		\bibitem{Opial1967}
		Z. Opial.
		\textit{Weak convergence of the sequence of successive approximations for nonexpansive mappings}.
		\emph{Bulletin of the American Mathematical Society}, \textbf{73}(4), 1967, pp. 591--597.


\bibitem{Park2022}
       J. Park, E. K. Ryu.
       \textit{Exact optimal accelerated complexity for fixed-point iterations}.
      \emph{ Proceedings of the 39th International Conference on Machine Learning,
       Proceedings of Machine Learning Research}, \textbf{162}, 2022, pp. 17420--17457.

   \bibitem{Chen2024}
        K. Chen, D. Sun, Y. Yuan, G. Zhang, X. Zhao.
       \textit{HPR-LP: An implementation of an HPR method for solving linear programming}.
      \emph{Mathematical Programming Computation}, 2025. https://doi.org/10.1007/s12532-025-00292-0.
       \bibitem{Lu2024}
        H. Lu, J. Yang.
    \textit{Restarted halpern PDHG for linear programming}.
    \emph{arXiv preprint}, arXiv:abs/2407.16144, 2024.

        \bibitem{Naldi2025}
		E. Naldi.
		\textit{Investigating degenerate preconditioners for proximal point algorithms}.
		\emph{Dissertationsschrift, TU Braunschweig, Braunschweig}, 2025.

		\bibitem{Xu}
		H. K. Xu.
		\textit{Iterative algorithms for nonlinear operators}.
		\emph{Journal of the London Mathematical Society}, \textbf{66}(1), 2002, pp. 240--256.
		
		\bibitem{Polyak}
		B. T. Polyak.
		\textit{Introduction to optimization}.
		\emph{Optimization Software, Inc.}, 1987.



		
		\bibitem{He2012}
		B. He, X. Yuan.
		\textit{Convergence analysis of primal-dual algorithms for a saddle-point problem: From contraction perspective}.
		\emph{SIAM Journal on Imaging Sciences}, \textbf{5}(1), 2012, pp. 119--149.
		
		\bibitem{Rockafellar1970}
		R. T. Rockafellar.
		\textit{Convex Analysis}.
		\emph{Princeton University Press, Princeton}, 1970.
		
		\bibitem{Chang2021}
		X. Chang, J. Yang.
		\textit{A golden ratio primal-dual algorithm for structured convex optimization}.
		\emph{Journal of Scientific Computing}, \textbf{87}(1), 2021, pp. 47--67.

		
		\bibitem{vonNeumann1928}
		J. von Neumann.
		\textit{Zur theorie der gesellschaftsspiele}.
		\emph{Mathematische Annalen}, \textbf{100}, 1928, pp. 295--320.

		
		\bibitem{lasso}
		R. Tibshirani.
		\textit{Regression shrinkage and selection via the lasso}.
		\emph{Journal of the Royal Statistical Society: Series B (Statistical Methodology)}, \textbf{58}(1), 1996, pp. 267--288.
		
		\bibitem{Duchi2008}
		J. Duchi, S. Shalev-Shwartz, Y. Singer, T. Chandra.
		\textit{Efficient projections onto the \(\ell_1\)-ball for learning in high dimensions}.
		\emph{Proceedings of the 25th International Conference on Machine Learning}, 2008, pp. 272--279.

        \bibitem{Ma2025}
       F. Ma, S. Li, and X. Zhang.
      \textit{Asymmetric version of the generalized Chambolle-Pock-He-Yuan method for saddle point problems}.
      \emph{Computational Optimization and Applications}, \textbf{91}, 2025, pp. 1--26.


	
	\end{thebibliography}
\end{document}